\documentclass[a4paper,12pt,reqno]{amsart}

\allowdisplaybreaks

\usepackage{rotating,pdflscape,color,amssymb,subfigure,psfrag,amsmath,eufrak,bbm,epsfig}
\usepackage{a4wide,amscd}
\usepackage{tikz} 
 \usepackage[usenames,dvipsnames]{pstricks}
 \usepackage{pst-grad} 
 \usepackage{pst-plot} 

\definecolor{vdarkred}{rgb}{0.6,0,0.2}
\definecolor{vdarkblue}{rgb}{0,0.2,0.6}
\usepackage[pdftex, colorlinks, linkcolor=vdarkblue,citecolor=vdarkred]{hyperref}
\usepackage[small]{caption}

\graphicspath{{../Figures/}}

 \DeclareMathOperator{\real}{Re}
\DeclareMathOperator{\im}{Im}
 \newcommand{\tM}{\wt{M}}
\newcommand{\dD}{\mathbb{D}}

\newcommand{\Cc}{\mc{C}}

\newcommand{\la}{\label}
 \newcommand{\eqre}{\eqref}
\newcommand{\re}{\ref}
\newcommand{\ld}{\ldots}
\newcommand{\beg}{\begin}
\newcommand{\en}{\end}

\newcommand{\trm}{\textrm}
 
\newcommand{\bgt}{\begin{itemize}}
\newcommand{\ent}{\end{itemize}}
\newcommand{\ite}{\item}

\newcommand{\op}{\operatorname}

\newcommand{\ds}{\displaystyle}

\newcommand{\p}{\mathbb{P}}

\newcommand{\supp}{\operatorname{supp}}

\newcommand{\Tr}{\operatorname{Tr}}

\newcommand{\Ninf}{\underset{N\to\infty}{\longrightarrow}}

\newcommand{\E}{\mathbb{E}}

\newcommand{\R}{\mathbb{R}}
\newcommand{\C}{\mathbb{C}}

\newcommand{\ud}{\mathrm{d}}

\newcommand{\pro}{probability }
\def\S{{\mathfrak{S}}}

\newcommand{\f}{\frac}
\newcommand{\ff}{\frac{1}}
\newcommand{\lf}{\left}
\newcommand{\ri}{\right}
 
\newcommand{\st}{such that }

\newcommand{\lam}{\lambda}

\newcommand{\ti}{\times}

\newcommand{\vfi}{\varphi}
\newcommand{\ste}{\, ;\, }
\newcommand{\mc}{\mathcal }

\newcommand{\eps}{\varepsilon}

\newcommand{\bxp}{\boxplus}

\newcommand{\A}{\mc{A}}

\newcommand{\B}{\mc{B}}

\newcommand{\al}{\alpha}
\newcommand{\tta}{\theta}

\newcommand{\eqlaw}{\stackrel{\textrm{law}}{=}}

\newcommand{\ovl}{\overline}

\newcommand{\bbm}{\begin{bmatrix}}
\newcommand{\ebm}{\end{bmatrix}}
\newcommand{\bes}{\begin{equation*}}
\newcommand{\ees}{\end{equation*}}
\newcommand{\be}{\begin{equation}}
\newcommand{\ee}{\end{equation}}
\newcommand{\beqy}{\begin{eqnarray}}
\newcommand{\eeqy}{\end{eqnarray}}
\newcommand{\beq}{\begin{eqnarray*}}
\newcommand{\eeq}{\end{eqnarray*}}
\newcommand{\one}{\mathbbm{1}}
\newcommand{\lto}{\longrightarrow}

\newcommand{\ie}{i.e. }

\newcommand{\bpm}{\begin{pmatrix}}
\newcommand{\epm}{\end{pmatrix}}

\newcommand{\wt}{\widetilde}

\newcommand{\bpr}{\beg{proof}}
\newcommand{\epr}{\en{proof}}

\newcommand{\bet}{\beta}
\newcommand{\del}{\delta}
\newcommand{\Del}{\Delta}

          \newcommand{\da}{\downarrow}

\newcommand{\U}{\mc{U}}

\newcommand{\pa}{\partial}

\newcommand{\bfa}{\mathbf{a}}
\newcommand{\bfb}{\mathbf{b}}

\newcommand{\bA}{\mathbf{A}}
\newcommand{\bB}{\mathbf{B}}
\newcommand{\bH}{\mathbf{H}}
\newcommand{\bW}{\mathbf{W}}

\newcommand{\bZ}{\mathbf{Z}}

 \newcommand{\ii}{\mathrm{i}}

\newcommand{\ka}{\kappa}

\newcommand{\oci}{\overset{\circ}}
\newcommand{\M}{\mc{M}_N(\C)}
\newcommand{\Md}{\mc{M}_{2N}(\C)}
 
  \newcommand{\mre}{\mathrm{e}}

\theoremstyle{definition}

\long\def\symbolfootnote[#1]#2{\begingroup
\def\thefootnote{\fnsymbol{footnote}}\footnote[#1]{#2}\endgroup}

\author[]{Florent Benaych-Georges} \address{MAP 5, UMR CNRS 8145 - Universit\'e Paris Descartes, 45 rue des Saints-P\`eres 75270 Paris cedex~6,  France.} \email{florent.benaych-georges@parisdescartes.fr}



\title[]{Local Single Ring Theorem}

\begin{document}
\maketitle

 \beg{abstract}The Single Ring Theorem, by Guionnet, Krishnapur and Zeitouni in \cite{GUI}, describes the   empirical eigenvalue   distribution   of a large generic matrix with prescribed singular values, \emph{i.e.} an $N\times N$  matrix of the form $A=UTV$, with $U, V$ some independent Haar-distributed unitary matrices and $T$ a deterministic matrix whose singular values are the ones prescribed. In this text, we give a \emph{local} version of this result,  proving that it remains true at the microscopic scale  $(\log N)^{-1/4}$. On our way to prove it, we prove a \emph{matrix subordination result} for singular values of sums of non-Hermitian  matrices, as Kargin did in \cite{KarginAOP13Sub}  for Hermitian matrices. This   allows to prove a local law for the singular values of the sum of two non-Hermitian matrices and a delocalization result for singular vectors.
 \end{abstract}
 
 

 \section*{Introduction}
The Single Ring Theorem, by Guionnet, Krishnapur and Zeitouni in \cite{GUI}, describes the   empirical eigenvalue  distribution  of a large generic matrix with prescribed singular values, \ie an $N\ti N$  matrix of the form $A=UTV$, with $U, V$ some independent Haar-distributed unitary matrices and $T$ a deterministic matrix whose singular values are the ones prescribed.    More precisely, under some technical hypotheses
, as the dimension $N$ tends to infinity, if  the  empirical distribution of the  singular values of $A$ (\ie of  $T$) converges to a compactly supported limit \pro measure $\nu$ on the real line, 
then the empirical eigenvalue  distribution of $A$ converges to a limit \pro measure $\mu$ on the complex plane
which depends only on $\nu$. The limit measure $\mu$ is radial, has support \be\la{bact15h20S}S:=\{z\in\C\ste a\le |z|\le b\}\quad\quad\trm{for }\; a^{-2}:=\int x^{-2}\ud\nu(x)\;;\,\; b^2:=\int x^{2}\ud\nu(x)\ee and  density $\rho$     satisfying \be\la{bact15h20}\rho (z):=\ff{2\pi}\Del_z(\int\log |x|\nu_{\infty,z}(\ud x)) \qquad\trm{with}\qquad \nu_{\infty,z}:= 
\nu^s\bxp\f{\del_{|z|}+\del_{-|z|}}{2}\ee  ($\nu^s$ is the symmetrization of $\nu$, see \eqre{941419h}, and $\bxp$ is the additive free convolution \cite{vdn91,ns06,agz}). 
\begin{figure}[ht]
\centering
\subfigure[Spectrum of the $500\times 500$  matrix $A=UTV$   when the singular values of $T$ are uniformly distributed on $(0 . 5 , 4)$]{
\includegraphics[scale=.35]{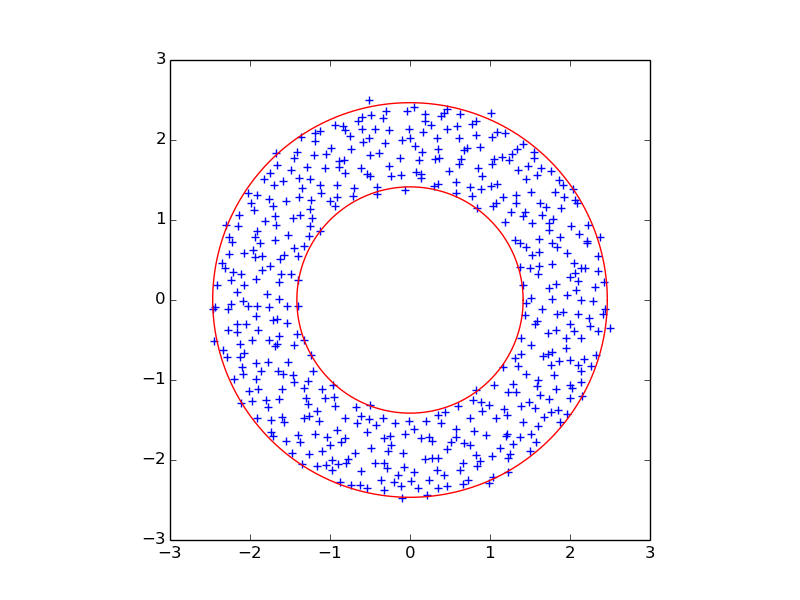} 
\label{FigCorNC1}} \qquad 
\subfigure[$500$ points uniformly distributed on the support of the limit spectral distribution of $A$]
{\includegraphics[scale=.35]{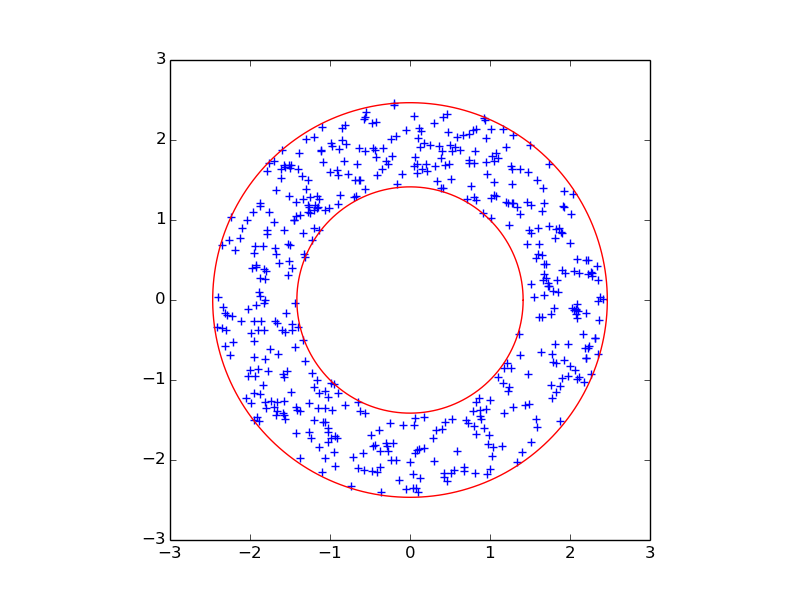} 
\label{Fig2CorC}}
\caption{{\bf Repulsion} (eigenvalues of $A$)/{\bf lack of repulsion} (independently distributed points)}\la{figure_intro}
\end{figure}
In the left image of Figure \re{figure_intro}, we plotted the spectrum of an example of such a matrix $A$ with size  $N=500$,  illustrating the convergence of the empirical spectral measure. In the right image of  Figure \re{figure_intro}, we plotted $500$ independent random points with uniform distribution on the ring $S$. Our point was not to compare the limit spectral distribution of $A$ with the uniform distribution on $S$, but to compare both point processes at microscopic scale: we see that the $500$ eigenvalues of $A$ fill the ring way more regularly    than the independent points, which reflects the so-called \emph{eigenvalues repulsion phenomenon}. Some of the mathematical manifestations of such phenomenons are the so-called \emph{local laws} (see e.g. \cite{ESY2,BourgadeCirc,BourgadeCirc2}). Here, we will prove a local law for the Single Ring Theorem on scale $\eps_N=(\log N)^{-1/4+\epsilon}$ in the interior of $S$, which means roughly that the number of eigenvalues of $A$ in any  ball $B(z_0,r)$ contained in $S$ is asymptotic to $\mu(B(z_0,r))\times N$ not only for fixed $r$   but also for $r\sim \eps_N$. 
 
 To give an idea of the techniques used in the proofs and of the difficulties we had to overcome, let us compare them with those of another  local law for non-Hermitian matrices. Recently, in    the series of papers \cite{BourgadeCirc,BourgadeCirc2,yinLCL3}, Bourgade, Yau and Yin proved a local law for non-Hermitian matrices with i.i.d. entries. It is well known that the empirical spectral distribution of a random matrix with size $N$ whose entries are i.i.d., centered, with variance $1/N$ and subject to no symmetry tends to the uniform measure on the unit disc of $\C$ when the dimension tends to infinity (see \cite{TaoVuKCL}). In    \cite{BourgadeCirc,BourgadeCirc2,yinLCL3}, the authors gave an almost optimal result about the local accuracy of the approximation of the empirical spectral distribution by its limit: they proved, through $\Cc^2$ test functions, that  the approximation stays correct as long as we consider test sets with surface  at least $N^{\epsilon-1}$, for any $\epsilon>0$. As a subset of the unit disc with normalized surface $\mc{S}$ should contain approximately $N\times \mc{S}$ eigenvalues, this is not far from the best one could do by considering sets with more than finitely many eigenvalues. 
 In the  local law we give here, we are far from this optimal scale, but the set of tools we have at disposal lacks several key elements. 
 The proofs, in  \cite{BourgadeCirc,BourgadeCirc2,yinLCL3} as well as  in the present paper,   are based on the  so-called \emph{Hermitization technique}, which expresses the empirical spectral distribution of a non-Hermitian  matrix $A$ as the Laplacian of the function $f(z)=\frac{1}{N}\Tr\log|A-z|$, with $|A-z|=\sqrt{(A-z)(A-z)^*}$ (see \eqre{251151}). 
 In \cite{BourgadeCirc,BourgadeCirc2,yinLCL3}, where $A$ is a matrix with i.i.d. entries,  $A-z$ is a matrix of the type ``information plus noise", a model well understood. It allows the authors   of  \cite{BourgadeCirc,BourgadeCirc2,yinLCL3} to prove, thanks to the Schur complements formula, that for any $z$, the empirical  eigenvalues distribution of $|A-z|$ is close to its limit at local scale $N^{\epsilon-1}$. 
 Then, as the limit spectral distribution of $|A-z|$ has a smooth density whose singularity points are well understood and as   the smallest singular values of $z-A$ are not likely to be too close to zero,
 the authors   of  \cite{BourgadeCirc,BourgadeCirc2,yinLCL3} approximate $\frac{1}{N}\Tr\log|A-z|$ by its theoretical limit quite well. 
 Here,    the Schur complements formula is not an option because suppressing a row and a column   breaks the symmetry of the Haar measure. Instead, we use  the \emph{matrix subordination}, a technique   proposed by    Kargin in  \cite{KarginAOP13Sub}: in  Theorem \re{241417h3NH}, we prove that for any matrix $B$ independent of $A=UTV$,  the resolvants  $G_\bA$, $G_\bB$ and $G_\bH$  of the matrices $$ \bA:=\bpm 0&A\\A^*&0\epm\qquad ;\qquad \bB:=\bpm 0&B\\B^*&0\epm\qquad \trm{ and }\qquad \bH:=\bpm 0&A+B\\(A+B)^*&0\epm$$ at $z=E+\ii\eta$, $\eta\ge N^{1/8}$, satisfy \beqy\la{2511513h} \E G_\bH(z)&=&G_\bA(z+S_B(z))+\lf(\trm{error term with operator norm $\le \ff{N\eta^6}$}\ri)\\ \la{2511513h1}
 \E G_\bH(z)&=&G_\bB(z+S_A(z))+\lf(\trm{error term with operator norm $\le \ff{N\eta^6}$}\ri)\eeqy for some complex-valued functions $S_{A},S_B$ \st $\im S_{A,B}(z)\ge -\ff{N\eta^7}$. Equations \eqre{2511513h} and \eqre{2511513h1} have to be compared with the ones   defining the free convolution $\bxp$ thanks to Stieltjes transforms subordination (see Theorem \re{24141}): \beqy\la{74141TD} m_{\mu^s\bxp\nu^s}(z)&=& m_{\nu^s}(z+S_{\mu}(z))\qquad; \qquad\im S_\nu(z)\ge 0,\\
 \la{74142TD}m_{\mu^s\bxp\nu^s}(z)&=& m_{\mu^s}(z+S_\nu(z))\qquad; \qquad\im S_\nu(z)\ge 0.\eeqy
The Hermitization technique described above brings us to use these equations with $B=-zI$ and  $\mu=\del_{|z|}$, so that  $\mu^s\bxp\nu^s=\nu_{\infty,z}$. 
Ideally, Equations \eqre{2511513h}--\eqre{74142TD} should give an upper-bound on $ \E \ff{2N}\Tr G_\bH(z)-m_{\mu^s\bxp\nu^s}(z)$ which could be turned into an upper-bound on \be\la{3011512h}\ff{2N}\Tr \log |\bH| -\int \log|x| \ud \mu^s\bxp\nu^s(x).\ee
The problem here is that   the upper-bound on $ \E \ff{N}\Tr G_\bH(z)-m_{\mu^s\bxp\nu^s}(z)$ deduced  from     Equations \eqre{2511513h}--\eqre{74142TD}   involves the inverse of a certain $2\ti 2$ determinant (see \eqre{14111423h}), which can vanish  for $z$ close to the real line (to control this determinant for $z$ close to the real line, one   would need precise informations about the density of $\nu_{\infty,z}$, which have, except for the case of i.i.d. matrices, remained out of reach so far, despite several studies of these questions as in \cite{BA,bbg07,bbcf12})\footnote{The lack of informations on the order of the density of $\nu_{\infty,z}$ at its singularities is also what makes the use of the estimates of Guionnet, Krishnapur and Zeitouni, like \cite[Eq. (12)]{GUI2}, ineffective here.}. However, using only  some bounds on the operator norms of $A$ and $B$, we can deduce from \eqre{2511513h}--\eqre{74142TD} that for $|z|$ large enough, \be\la{291152}\lf| \E \ff{N}\Tr G_\bH(z)-m_{\mu^s\bxp\nu^s}(z)\ri|\le \f{C}{N}.\ee The necessity to  have $|z|$ large   for such a bound to be proved is a real problem in the perspective of  establishing a local law for the eigenvalues of $A$. We   fix it (at the price of a quite poor microscopic scale $\eps_N$) 
  using Hadamard's three circles theorem, an idea introduced by    Kargin in \cite{KarginPTRF12}. This theorem, with some standard concentration inequalities, allows to deduce from \eqre{291152} that for $\eta \sim \ff{\sqrt{\log N}}$, we have \be\la{291151}   \ff{N}\Tr G_\bH(z)-m_{\mu^s\bxp\nu^s}(z) \ll 1.\ee  
To conclude the proof, we need to turn \eqre{291151} into a control on \eqre{3011512h}: this is done thanks to the    Helffer-Sj\H{o}strand functional calculus  and to a recent theorem by Rudelson and  Vershynin in  \cite{RUD} about the smallest singular value of $A-z$.\\
The recent preprints \cite{BaoErdosSchnelli1,BaoErdosSchnelli2} by Bao, Erd\H{o}s and Schnelli give local laws for the close model $A+UBU^*$ when $\mu_A\to\mu_\al$ and $\mu_B\to\mu_\bet$ as $N\to\infty$. Their local laws are established at some better scales than the ones we give here for $UTV^*-z$, but, seemingly facing the same problem as us, they had to specify the part $\mc{B}_{\mu_\al\bxp\mu_\bet}$ of the real line where they establish these laws, avoiding a set of singular points (see  \cite[Th. 2.7]{BaoErdosSchnelli1} and \cite[Th. 2.5]{BaoErdosSchnelli2}). It should be possible to adapt their proofs to our  model $UTV^*-z$, but at the current level of understanding of the densities of the  laws $\nu_{\infty,z}$, we do not know exactly what the sets $\mc{B}_{\nu_{\infty,z}}$ look like and how to deal with their complementaries, thus it is today not possible to convert such local laws into a local version of the Single Ring Theorem. 
\\

\noindent{\bf Organization of the article :} We postpone the proof of our key result, the matrix subordination result (Theorem \re{241417h3NH})  to Section \re{PrMaSub200115}. We will first prove its main consequence, Proposition \re{18111411h1}, in Section \re{20115.18h}. Then, the short proofs of the  local law for the singular values of $A+\wt{B}$ (Theorem \re{241417h5}) and of the singular vectors delocalization for  $A+\wt{B}$ (Theorem \re{241417h6}) will be given in Section \re{prfreeresults}. The proof of the local Single Ring Theorem will be given in Section \re{secPrLSRT}, followed in Section \re{PrMaSub200115} by the proof of Theorem \re{241417h3NH} and in the appendix by several results we will use here.
 
 \noindent {\bf Notation :} Throughout this text, $z=E+\ii\eta$, $E\in \R$, $\eta>0$ denotes an element of $\C^+:=\{\xi\in \C\ste \im \xi>0\}$.   For $\mu$ a  signed  measure  on the real line, we define \be\la{941419h}\mu^{s}(X):=\ff{2}(\mu(X)+\mu(-X))\qquad\trm{(\emph{symmetrization} of $\mu$)},\ee \be\la{941419h2}m_\mu(z):=\int \f{\ud\mu(t)}{t-z}\qquad \trm{(\emph{Stieltjes transform} of $\mu$)}\ee and for $M\in \M$ (the set of $N\ti N$ complex matrices), we define $\mu_M$ as the \emph{empirical eigenvalue  distribution} of $M$: \be\la{941419h3}\mu_M:=\ff{N}\sum_{i=1}^{N} \del_{\lam_i}\qquad\trm{($\lam_1, \ld, \lam_{N}$: eigenvalues of $M$)},\ee
 whereas $\nu_M$ denotes the \emph{empirical singular value    distribution} of $M$:   \be\la{941419h5}\nu_M:=\ff{N}\sum_{i=1}^N \del_{s_i}\qquad\trm{($s_1, \ld, s_N$: singular values of $M$)}.\ee  Note that we have \be\la{241417hMLV012}\nu^{s}_M=\mu_{\mathbf{M}}\quad\quad \trm{ for }\quad\mathbf{M}:=\bpm 0&M\\ M^*&0\epm.\ee
 We   denote by $\|M\|$ the canonical operator norm of $M$.
 When $M$ is Hermitian, we also define, for  $z=E+\ii\eta,\;E\in \R,\;\eta>0$, \be\la{941419h4} m_M(z):=m_{\mu_M}(z)=\ff{N}\Tr G_M(z) , \ee for  $G_M$ the \emph{resolvant matrix} of $M$:  \be\la{201141}G_M(z)=(M-z)^{-1}  .\ee
 

 \noindent For $X$ an $L^1$ random variable,  \be\la{941416h32}\oci{X}:=X-\E X.\ee Note that for any $X,Y\in L^2$, \be\la{941416h33}\E[XY]=\E X \E Y+\E[\oci{X}Y].\ee
 
 \noindent For $f$ a function of a real variable and $\ell\ge 0$, $f^{(\ell)}$ denotes the $\ell$-th derivative of $f$.
 
 \noindent For $E\in \R$ and $\del>0$, $[E\pm\delta]$ denotes the interval $[E-\del, E+\del]$.
 
 \noindent For $X=X_N$ and  $Y=Y_N$, $X\ll Y$ means that $X/Y\lto 0$ as $N\to\infty$.

 \section{Main results}
 \subsection{Local Law for the Single Ring Theorem}Let $A$ be an $N\ti N$  matrix (depending implicitly on $N$) of the form $A=UTV$, with $U, V$ some independent Haar-distributed unitary matrices and $T$ a deterministic matrix. 
   
   We make the following hypothesis. 

\beg{hyp}\la{2711141}\bgt\ite[(i)] There is $K>0$ independent of $N$ \st $\|T\|\le K$,
\ite[(ii)] There is $\nu$ a \pro measure, $\eta_0>0$ and    $C_0$ independent of $N$  \st $$ \eta=\im z>\eta_0\implies |m_{\nu_T-\nu}(z)|\le C_0N^{-1},$$ where $m$ is the Stieltjes transform defined at \eqre{941419h2} and  \eqre{941419h4},
\ite[(iii)] There are $C,c>0$ independent of $N$ \st $\im m_{\nu_T}(z)\le C$ when $ \im z>N^{-c}$.\ent
 \en{hyp}
 
 Then we know, by \cite{GUI,RUD}\footnote{For the Single Ring Theorem to hold, these hypotheses can even be weakened, as proved  by Basak and  Dembo in \cite{BD13}.}, that $\mu_A$ converges in \pro to a law $\mu$ with density $\rho$ given by \eqre{bact15h20} and support $S=\{z\in \C\ste a\le |z|\le b\}$ given by \eqre{bact15h20S}.

    Here is our main result. It will be proved in Section \re{secPrLSRT} as a consequence   of the local law for the singular values proved in Theorem \re{132151}, which is in turn proved   using the matrix subordination result in Theorem \re{241417h3NH} and its consequence in Proposition \re{18111411h1}.

\beg{Th}[Local Single Ring Theorem]\la{SRT}
Fix  $z_0$ \st $a<|z_0|<b$, $\al\in (0,1/4)$ and define $\eps_N:=(\log N)^{-\al}$. Then     for  $f\in \Cc^2_c(\C)$ and   $$F_{z_0,\eps_N}:\lam\longmapsto \eps_N^{-2}f( \f{\lam-z_0}{\eps_N}),$$ we have the convergence in probability, as $N\to\infty$, $$\int F_{z_0,\eps_N}(\lam)\ud \mu_A(\lam)-\int F_{z_0,\eps_N}(\lam)\ud\mu(\lam)\lto 0,$$ where $\mu$ is the limit spectral law of $A$, introduced above.
\en{Th} 
 
 \beg{rmk}Why do  we call it a   \emph{local law}?
The convergence of $\mu_A$ towards  $\mu$  can be considered as  local  with scale $\eps_N$   at $z_0
$ when for any    test function $f$, \be\la{WDWCIALL}   \int f( \f{\lam-z_0}{\eps_N})\ud \mu_A(\lam)- \int f( \f{\lam-z_0}{\eps_N})\ud\mu(\lam)\;\ll\;   \int f( \f{\lam-z_0}{\eps_N})\ud\mu(\lam).\ee As for a test function $f$ with enough decay at infinity,  the RHT of \eqre{WDWCIALL} should have order at most   $ \mu (B(z_0,  \eps_N))\approx \eps_N^2$, this rewrites   $$\int f( \f{\lam-z_0}{\eps_N})\ud \mu_A(\lam)- \int f( \f{\lam-z_0}{\eps_N})\ud\mu(\lam)\;\ll\; \eps_N^2,$$ which is precisely the contents of the theorem.\en{rmk}
 \beg{rmk} Note that we focus on the interior of the support $S$ (it is necessary at \eqre{2611513h}). It has been proved in \cite{GUI2,BEN} that there is no eigenvalue at a macroscopic distance of $S$, but the case of the border of $S$ (\ie$|z_0|=a$ or $b$) is not treated here.
 \end{rmk}
 
 
  \subsection{Matrix subordination}
 In order to prove this local law, we need to prove a \emph{matrix subordination result}, as  called by Kargin in \cite{KarginAOP13Sub}, where he introduced this idea for Hermitian matrices. 
 Let $A,B$  be deterministic, depending implicitly on $N$, $N\ti N$ matrices    \st there is $K$ independent of $N$ \st 
\be\la{941416h112}\|A\|, \|B\|\le K.\ee Let $U, V$ be some independent $N\ti N$ Haar-distributed unitary matrices and \be\la{164146h}\wt{B}:=UBV^*.\ee
We introduce the matrices  \be\la{1141413h}\bA:=\bpm 0&A\\ A^*&0\epm\qquad\bB:=\bpm 0&B\\ B^*&0\epm\qquad\bW:=\bpm U& 0\\ 0& V\epm\qquad \wt{\bB}:=\bW\bB\bW^*\ee \be\la{1141413h1}\bH:=\bpm 0&A+\wt{B}\\ (A+\wt{B})^*&0\epm   =\bA+\wt{\bB},\ee Note that the matrices $\bA$, $\bB$ and $\bH$ have eigenvalues the singular values of respectively $A$, $B$ and  $A+\wt{B}$ (and their opposites).

  \beg{Th}\la{241417h3NH}There are some complex valued functions $S_A(z)$, $S_B(z)$ of $z=E+\ii\eta\in \C^+$  and some matrices $R_A(z)$ and $R_B(z)$ \st we have \beqy\la{34143} \E G_\bH(z)&=&G_\bA(z+S_B(z))+R_A(z)\\ \nonumber \\ \la{34144}\E G_\bH(z)&=&G_\bB(z+S_A(z))+R_B(z),\eeqy
 \st  the  functions $S_A(z)$, $S_B(z)$, whose formulas are given at \eqre{941416h12}, satisfy
  \be\la{241417h}N\eta^5\ge C\implies \im S_A(z), \im S_B(z)\ge -\f{C }{N\eta^7}\ee and \st  the matrices  $R_A(z)$, $R_B(z)$,  whose formulas are given at \eqre{941416h12irmar}, satisfy
   \be\la{241417h2}N\eta^8\ge C \implies \|R_A(z)\|, \|R_B(z)\|\le \f{C }{N\eta^6}\ee for a constant $C$ depending only on the $K$ of \eqre{941416h112}.
  \en{Th}
  
  \beg{rmk}Theorem \re{241417h3NH}, which will be proved in Section \re{PrMaSub200115},  has to be compared with Theorem \re{24141} and Remark \re{201152} of the appendix, which give the definition of the free convolution $\bxp$ in terms of \emph{subordination of Stieltjes transforms} and interpret subordination as a \emph{regularity criterion}. 
  \en{rmk}
  
  \subsection{Local laws for the singular values of $A+B$ and singular vectors delocalization}
On the way to deduce Theorem \re{SRT} from Theorem \re{241417h3NH}, we will prove a key result, Proposition \re{18111411h1}. Then, for almost free, we get the two following results  (Theorems \re{241417h5}  and \re{241417h6}). 

Let us suppose that besides the  hypothesis $\|A\|, \|B\|\le K$, there are some \pro measures $\nu_\bfa\ne\del_0$ and $\nu_{\bfb}\ne\del_0$ \st as $N\to\infty$,  
 \be\la{1911416h41} \nu_A\lto \nu_{\bfa}\qquad ;\qquad \nu_B\lto \nu_{\bfb}.\ee  
It is well known \cite{haag2,BENRECT} that we then have the convergence in \pro  \be\la{241417hMLV}\nu_{A+\wt{B}}\lto \nu,\ee  with $\nu$ the \pro measure on $\R_+$ whose symmetrization $\nu^{s}$ (see \eqre{941419h}) 
satisfies  \be\la{164146h2}\nu^{s}=\nu_{\mathbf{a}}^{s}\bxp\nu_{\mathbf{b}}^{s}.\ee
The two next theorems give conditions for    the convergence of \eqre{241417hMLV} to hold at local levels. 

\beg{Th}[Local law 1 for the singular values of $A+\wt{B}$]\la{132151}Suppose that there are $\eta_0,C_0$ independent of $N$ \st $$\eta\ge \eta_0\implies |m_{\nu_A-\nu_{\bfa}}(z)|+|m_{\nu_B-\nu_{\bfb}}(z)|\le C_0 N^{-1}.$$ Let $p\ge 0$ and let $\phi_N$ be a sequence of smooth      functions.  Then there are  $C  , c>0$ depending only on  $K,\nu_{\bfa},\nu_{\bfb},C_0,\eta_0,p$ \st  with \pro at least $1-C \mre^{-N^c}$, $$|(\nu_{A+\wt{B}}^s-\nu_{\mathbf{a}}^{s}\bxp\nu_{\mathbf{b}}^{s})(\phi_N)|\le \f{C \|\phi_N^{(p+1)}\|_\infty}{(\log N)^{p/2}}.$$
\en{Th}

 We define, for $\mu,\nu$, compactly supported \pro measures on $\R$   and $z\in \C^+$, \beqy \nonumber\ka_{\mu,\nu}(z)&:=&\{m_{\mu}'(z+S_\nu(z))+m_{\nu}'(z+S_\mu(z))\}(z+S_\mu(z)+S_\nu(z))^{-2}\\&&\qquad\qquad\qquad\qquad\qquad\qquad -m_{\mu}'(z+S_\nu(z))m_{\nu}'(z+S_\mu(z)),\la{14111423h} \eeqy
where the functions $S_\mu$, $S_\nu$ are the subordination functions introduced in Theorem \re{24141}.

 We use the definition introduced by Kargin in \cite{KarginAOP13Sub}: \beg{Def}\la{211141}We say that the pair $(\mu,\nu)$ of \pro measures on $\R$  is \emph{well behaved} at $E\in \R$ if:\bgt\ite[a)] the  subordination  functions $S_\mu$, $S_\nu$ have finite limits\footnote{It has been proved in \cite{SerbanLinfty} that for $\mu,\nu$ compactly supported, $S_\mu$, $S_\nu$ extend continuously to the whole real line, with values in $\C^+\cup\R\cup\{\infty\}$.} with positive imaginary parts at $E$,\\
 \ite[b)] the value of the analytic continuation\footnote{It has been proved at Th. 3.3 of \cite{BA} that a) implies that the functions $S_\mu$ and $S_\nu$ have analytic continuations to a neighborhood of $E$, which implies   that $\ka_{\mu,\nu}(z)$ does so.} of the function $\ka_{\mu,\nu}(z)$ at $E$ is non zero.\\
 \ent\end{Def}
 
 \beg{rmk}\la{211152h52}Sufficient conditions have been given, for example by Belinschi in \cite{BA}, for a) of the previous definition to occur. As far as b) is concerned, Kargin gave sufficient conditions in \cite{KarginAOP13}. Besides, if a) is satisfied, by the analytic continuation principle and an analysis of the function $\ka_{\mu,\nu}(z)$ at infinity, we see that the set of $E$'s where a) holds and not b) is discrete. 
 \en{rmk}

  \beg{Th}[Local law 2 for the singular values of $A+\wt{B}$]\la{241417h5} Suppose that  the pair $(\nu_\bfa^s,\nu^s_\bfb)$ is well behaved  at $E\in \R$ and that there is $\eta_0=\eta_0(N)$ \st $$\sup_{\eta\ge \eta_0} \eta^2|m_{\nu_A-\nu_{\bfa}}(z)|+\eta^2|m_{\nu_B-\nu_{\bfb}}|\;\ll \;   N^{-1/4}.$$
 Then  we have      $$\max\{\eta_0,N^{-1/8}\}\ll \eta\ll 1\implies \f{\nu_{A+\wt{B}}([E\pm\eta])}{2\eta}\lto \rho(E)$$ for the convergence in probability, where $\rho$ is the density\footnote{It follows from Th. 7.4 of \cite{BV1} and  Th. 4.1 p. 146 of \cite{BA} that there is an open set $U\subset \R^+$ and an analytic positive function $\rho$ on $U$ \st the limit $\nu$  of $\nu_{A+\wt{B}}$ is $  \al\del_0+\one_U(x)\rho(x)\ud x$ for $\al:=((\nu_\bfa+\nu_\bfb)(\{0\})-1)_+$.}, at $E$, of the limit of $\nu_{A+\wt{B}}$, \ie   of  $\nu_\bfa^s\bxp\nu^s_\bfb$.
 \en{Th}

 \beg{rmk}The statement of Theorem \re{241417h5} is close, in nature, from the one of Theorem \re{132151}. However, the statement of Theorem \re{241417h5} (local law at scale $N^{-1/8}$) is  stronger than the one of  Theorem \re{132151} (local law at logarithmic scale), but relies on stronger hypotheses (we need to know that $(\nu_\bfa^s,\nu^s_\bfb)$ is well behaved  at $E\in \R$, which is usually hard to prove, given how little explicit formulas for $\bxp$ are). This dichotomy 
 is reflected in an essential difference in their proofs: the proof of Theorem \re{241417h5} relies on Erd\H{o}s, Schlein and Yau's method via the approximation of the Stieltjes transform of $\nu_\bfa^s\bxp\nu^s_\bfb$ by the one of $\nu_{A+\wt{B}}$ at distance $\eta$ from the real line (see Theorem \re{methodErdosSchleinYau2009}), whereas the proof  of Theorem \re{132151} relies on Hadamard's three circles theorem and the approximation of the Stieltjes transform of $\nu_\bfa^s\bxp\nu^s_\bfb$ by the one of $\nu_{A+\wt{B}}$ at macroscopic distance  from the real line (see Corollary \re{corH3CTh}).
 \en{rmk}
 
 Let us now state a result about the delocalization of the singular vectors of $A+\wt{B}$ which will also come for free once 
 Theorem \re{241417h3NH} and Proposition \re{18111411h1} proved.
 
 Let
 $s_a$ ($a=1, \ld, N$) denote the singular values of $A+\wt{B}$ and let $u_a$  ($a=1, \ld, N$), $v_a$  ($a=1, \ld, N$)  denote some orthonormal bases \st for each $a$, $(u_a, v_a)$ is a pair of singular vectors for $A+\wt{B}$ associated to $s_a$ (\ie $A+\wt{B}=\sum_{a} s_au_av_a^*$). For each $a,i$, $u_a(i)$, $v_a(i)$ denote the $i$-th components of $u_a$, $v_a$. 
  \beg{Th}[Singular vectors delocalization for  $A+\wt{B}$]\la{241417h6}If the pair $(\nu_\bfa^s,\nu^s_\bfb)$ is well behaved  at each point in $[E-\eps,E+\eps]$ ($E\in \R$, $\eps>0$) and the hypotheses of Theorem \re{241417h5} hold, then  we have   
  $$ \p(\exists a,i \ste |\lam_a-E|\le \eps\trm{ and }( |u_a(i)|^2>CN^{-1/8}\log(N)\trm{ or }|v_a(i)|^2>CN^{-1/8}\log(N)))$$ $$\ \le \ \mre^{-c\sqrt{N}},$$ for some constants $c,C$ depending only on the parameters of the hypotheses.  
 \en{Th}

 \noindent{\bf Note about the constants $c,C$ :} In the proof of the  Local Single Ring Theorem (Theorem \re{SRT}) $c,C$ will denote some respectively small and large constants that might change from line to line and that depend only on the constant parameters introduced in the statement of Theorem \re{SRT} and in  Hypothesis \re{2711141}. In the same way, in the proofs of the matrix subordination result (Theorem \re{241417h3NH}),  the local law for singular values and the singular vectors delocalization (Theorems \re{241417h5} and  \re{241417h6}), as well as    Proposition \re{18111411h1}, $c,C$ might change from line to line and depend only on the parameters introduced in the hypotheses.

 \section{Statement and proof of Proposition \re{18111411h1}}\la{20115.18h}

For $\mu,\nu$, compactly supported \pro measures on $\R$   and $z\in \C^+$, besides the number $\ka_{\mu,\nu}(z)$ defined at \eqre{14111423h}, when $\ka_{\mu,\nu}(z)\ne 0$, we define the numbers \beqy 
 \la{defalpha}\al_{\mu,\nu}(z)&:=&\f{|z+S_{\mu}(z)+S_{\nu}(z)|^{-2}+|m_{\mu}'(z+S_{\nu}(z))|+|m_{\nu}'(z+S_{\mu}(z))|}{|\ka_{\mu,\nu}(z)|}\\ \nonumber&&\\
\la{defbeta} \bet_{\mu,\nu}(z)&:=& |z+S_{\mu}(z)+S_{\nu}(z)|^{-3}+|m_{\mu}''(z+S_{\nu}(z))|+|m_{\nu}''(z+S_{\mu}(z))|,\eeqy
where the functions $S_\mu$, $S_\nu$ are the subordination functions introduced in Theorem \re{24141}.

  The following consequence of Theorem \re{241417h3NH} will be a key result in the proof of the local version of the single ring theorem. Kargin  stated very similar results in \cite{KarginAOP13,KarginAOP13Sub} but to prove the local Single Ring Theorem, we need to give more accurate upper bounds than the ones given in Kargin's works. We sill use this proposition for small vaues of  $\eta$ in the proofs of Theorems \re{241417h5} and \re{241417h6} and macroscopic, as large as needed, values of $\eta$ in the proof of Theorem \re{132151}, which is a key step in the proof of the Local Single Ring Theorem.
     \beg{propo}\la{18111411h1}     Let $s\in (0,c)$   and $z=E+\ii\eta\in \C^+$ be \st   $$\ka_{\nu_{\mathbf{a}}^s,\nu_{\mathbf{b}}^s}(z)\ne 0\quad;\quad
     N\eta^8\ge c\quad;\quad N\eta^{6}\ge\f{c}{\al_{\nu_{\mathbf{a}}^s,\nu_{\mathbf{b}}^s}(z)\bet_{\nu_{\mathbf{a}}^s,\nu_{\mathbf{b}}^s}(z)},$$ 
     $$\forall z'\in \C, \qquad \im z'\ge \eta\implies | m_{\nu_{\mathbf{a}}}(z')-m_{ \nu_A}(z')|\le s, $$
     then  
   the following inequalities hold: \beqy\la{1411141} \lf|S_A(z)-S_{\nu_{\mathbf{a}}^s}(z)\ri|&\le &c\, \al_{\nu_{\mathbf{a}}^s,\nu_{\mathbf{b}}^s}(z)(  (N\eta^6)^{-1}+ s)\\
   \la{1411142} \lf|S_B(z)-S_{\nu_{\mathbf{b}}^s}(z)\ri|&\le &c\, \al_{\nu_{\mathbf{a}}^s,\nu_{\mathbf{b}}^s}(z)(  (N\eta^6)^{-1}+ s)\\
  \la{1411143}\lf|\E m_{\bH}(z)-m_{\nu_{\mathbf{a}}^s\bxp\nu_{\mathbf{b}}^s}(z)\ri|&\le& \f{c\, \al_{\nu_{\mathbf{a}}^s,\nu_{\mathbf{b}}^s}(z)(  (N\eta^6)^{-1}+ s)}{|z+S_{\nu_{\mathbf{a}}^s}(z)+S_{\nu_{\mathbf{b}}^s}(z)|^2}  .\eeqy
   \end{propo}

  \bpr Note first that, by  Theorem \re{241417h3NH} (whose proof is postponed to Section \re{PrMaSub200115}), \be\la{131411Andu}\E m_{\bH}(z)=\ff{2N}\Tr(\E G_\bH(z))=-\ff{z+S_A(z)+S_B(z)}.\ee Indeed,   with our definition of \eqre{941416h12} (and its analogue for $S_A(z)$), we have \beq z+S_A(z)+S_B(z)&=&\f{\E[zm_\bH(z)-f_B(z)-f_A(z)]}{\E[m_\bH(z)]}\\
&=&\f{\E[\ff{2N}\Tr((z-\bA-\bW\bB\bW^*)(\bA+\bW\bB\bW^*-z)^{-1})]}{\E[m_\bH(z)]}\\
&=&-\ff{\E[m_\bH(z)]}, 
\eeq
Besides, $$m_{\nu_{\mathbf{a}}^s\bxp\nu_{\mathbf{b}}^s}(z)=-\ff{z+S_{\nu_{\mathbf{a}}^s}(z)+S_{\nu_{\mathbf{b}}^s}(z)},$$
hence \beq\E m_{\bH}(z)-m_{\nu_{\mathbf{a}}^s\bxp\nu_{\mathbf{b}}^s}(z)&=&\f{S_{\nu_{\mathbf{a}}^s}(z)+S_{\nu_{\mathbf{b}}^s}(z)-S_A(z)-S_B(z)}{(z+S_A(z)+S_B(z))(z+S_{\nu_{\mathbf{a}}^s}(z)+S_{\nu_{\mathbf{b}}^s}(z))} ,\eeq
and the third equation of the lemma follows from its two first ones. Let us prove them.

   For  $z\in \C^+$, we define the set $$O_z:=\{(s_1,s_2)\in \C\ste s_1+s_2\ne -z\,,\;  z+s_1, z+s_2\in \C^+\}$$ and for $\mu,\nu$ \pro measures on $\R$, we define  the function $F_{\mu,\nu,z}:O_z\to \C^2$   by $$  F_{\mu,\nu,z}\bpm s_1\\s_2\epm:=\bpm m_{\mu^s}(z+s_2)+(z+s_1+s_2)^{-1}\\ m_{\nu^s}(z+s_2)+(z+s_1+s_2)^{-1}\epm.$$ 
  With the notations of  Theorem \re{24141}, we have \be\la{1311141}F_{\mu,\nu,z}\bpm S_{\mu^s}(z)\\ S_{\nu^s}(z)\epm=\bpm 0\\ 0\epm.\ee
  We shall apply it for $\mu=\nu_{\mathbf{a}}$ and $\nu=\nu_{\mathbf{b}}$, yielding \be\la{1311141bis}F_{\nu_{\mathbf{a}},\nu_{\mathbf{b}},z}\bpm S_{\nu_{\mathbf{a}}^s}(z)\\ S_{\nu_{\mathbf{b}}^s}(z)\epm=\bpm 0\\ 0\epm.\ee
    A similar system can be written for $S_A$ and $S_B$:  by \eqre{131411Andu},  \eqre{34143} and \eqre{34144} \beqy\la{34143J} m_\bA(z+S_B(z))+\ff{ z+S_A(z)+S_B(z)}&=&-\ff{2N}\Tr R_A(z) \ =: \ r_A(z)\\ \nonumber \\ \la{34144J}m_\bB(z+S_A(z))+\ff{ z+S_A(z)+S_B(z)}&=&-\ff{N}\Tr R_B(z)\ =: \ r_B(z),\eeqy
so that \be\la{1311142ante}F_{\nu_A,\nu_B,z}\bpm S_A(z)\\ S_B(z)\epm=\bpm r_A(z)\\ r_B(z)\epm\ee and that by hypothesis, 
\be\la{1311142}F_{\nu_\bfa,\nu_\bfb,z}\bpm S_A(z)\\ S_B(z)\epm=\bpm \hat{r}_A(z)\\ \hat{r}_B(z)\epm
  \trm{ 
with $|\hat{r}_A(z)-r_A(z)|+|\hat{r}_B(z)-r_B(z)|\le   s$.}
\ee

Let us now consider the  intermediate system:\be\la{1311143}F_{\nu_{\mathbf{a}},\nu_{\mathbf{b}},z}\bpm \tilde{S}_A(z)\\ \tilde{S}_B(z)\epm=\bpm r_A(z)\\ r_B(z)\epm.\ee
 Firstly, we shall  upper-bound the distance between the solution $(S_{\nu_{\mathbf{a}}^s}(z), S_{\nu_{\mathbf{b}}^s}(z))$ of \eqre{1311141bis} and the solution $(\tilde{S}_A(z), \tilde{S}_B(z))$  of \eqre{1311143}  (using     Kantorovich's Theorem \re{kanto}  and the fact that the derivative of $F_{\nu_A,\nu_B,z}$ is not too small and that $r_A(z), r_B(z)$ are small). Secondly, we shall    upper-bound the distance between  $(\tilde{S}_A(z), \tilde{S}_B(z))$  (as a solution of    \eqre{1311143} again) and the solution  $( {S}_A(z),  {S}_B(z))$ of \eqre{1311142}  (using the same ideas). 

Let us upper-bound the distance between the solutions of \eqre{1311141bis} and \eqre{1311143} thanks to Kantorovich's Theorem \re{kanto}. For $$S_z:=\bpm S_{\nu_{\mathbf{a}}^s}(z)\\ S_{\nu_{\mathbf{b}}^s}(z)\epm\qquad;\qquad M_z:=\|(F'_{\nu_{\mathbf{a}},\nu_{\mathbf{b}},z}(S_z))^{-1}\|,$$ 
 by Theorem \re{kanto}, we have $$|S_A(z)-S_{\nu_{\mathbf{a}}^s}(z)|+|S_B(z)-S_{\nu_{\mathbf{b}}^s}(z)|\le 100 M_z (|r_A(z)|+|r_B(z)|)$$ as soon as $$100M_z^2 (|r_A(z)|+|r_B(z)|) \|F''_{\nu_A,\nu_B,z}(S_z)\|<1$$   (the $100$ is here to avoid any norm choice issue, as  Theorem \re{kanto} is stated for the euclidian norm).
The derivative $$F_{\nu_{\mathbf{a}},\nu_{\mathbf{b}},z}'\bpm s_1\\s_2\epm=\bpm-(z+s_1+s_2)^{-2}&-(z+s_1+s_2)^{-2}+m_{\nu^s_{\mathbf{a}}}'(z+s_2)\\ -(z+s_1+s_2)^{-2}+m_{\nu^s_{\mathbf{b}}}'(z+s_1)&-(z+s_1+s_2)^{-2}\epm$$has determinant 
$$ \det F_{\nu_{\mathbf{a}},\nu_{\mathbf{b}},z}'\bpm s_1\\s_2\epm=\{m_{\nu_{\mathbf{a}}^s}'(z+s_2)+m_{\nu_{\mathbf{b}}^s}'(z+s_1)\}(z+s_1+s_2)^{-2}-m_{\nu_{\mathbf{a}}^s}'(z+s_2)m_{\nu_{\mathbf{b}}^s}'(z+s_1), $$
so that $$M_z\le \al_{\nu_{\mathbf{a}}^s,\nu_{\mathbf{b}}^s}(z)$$
Its second derivative is bounded by   $$100\bet_{\nu_{\mathbf{a}}^s,\nu_{\mathbf{b}}^s}(z).$$
This proves that under the hypotheses of the lemma, the distance between the solutions of \eqre{1311141bis} and \eqre{1311143}, \ie between $\bpm S_{\nu_{\bfa}^s}(z)\\S_{\nu_{\bfb}^s}(z)\epm$ and $\bpm \tilde{S}_{A}(z)\\ \tilde{S}_{B}(z)\epm$, 
is upper-bounded by the first part of the common right hand side of   \eqre{1411141} and \eqre{1411142}. 

Upper-bounding the distance between the solutions of \eqre{1311143} and \eqre{1311142} goes along the same lines, and gives the second part (the one involving $s$) of the common right hand side of   \eqre{1411141} and \eqre{1411142}. 
\epr

%

 \section{Proofs of Theorems  \re{132151}, \re{241417h5}  and \re{241417h6} (Local laws and singular vectors delocalization for $A+\wt{B}$)}\la{prfreeresults}
 
  \subsection{Proof of Theorem \re{132151} (Local law 1)}
  \beg{lem}\la{18111411h}Let   $\mc{K}$ be a fixed compact subset of $\C^+$.  Then there is $C=C(\mc{K})>0$ \st for any $\del>0$,     $$\p(\sup_{z\in \mc{K}}|m_{\bH}(z)-\E m_{\bH}(z)|>\del)\le C\mre^{-C \del^2N^2}$$
\end{lem} \bpr The lemma can be proved as Corollary 6 of \cite{KarginPTRF12}.\epr
  
  Let us now prove Theorem \re{132151}. 
  By Lemma \re{231151}, Proposition \re{18111411h1}, Lemma \re{18111411h} and  Corollary \re{corH3CTh}, we know that there are $C ,c>0$ \st    we have,  with  probability at least $1-C \mre^{-N^c}$,  $$\sup  \{|m_{\nu_{A+\wt{B}}^s-\nu_{\mathbf{a}}^{s}\bxp\nu_{\mathbf{b}}^{s}}(z)|\ste z=E+\ii\eta,\, |E|\le 3K,\, \eta= \f{C}{\sqrt{\log N}}\} \le C\mre^{-c\sqrt{\log N } }.$$
By   Corollary \re{2611141}, we conclude.   
   \subsection{Proof of Theorem \re{241417h5}  (Local law 2)}
   The proof  is a direct application of Lemma \re{18111411h} and  Proposition \re{18111411h1}   and of  Erd\H{o}s, Schlein and Yau's method (see Theorem \re{methodErdosSchleinYau2009} in the appendix).

  \subsection{Proof of Theorem \re{241417h6} (Singular vectors delocalization)}The proof is a copy of the one of Theorem 4 in  \cite{KarginAOP13}. Let us give the main lines. 
 First, we note that for any  $a$, as $(u_a, v_a)$ of unit singular vectors associated to the singular value $s_a$  of $A+\wt{B}$, the vector $$w_{a, \pm}:=\ff{\sqrt{2}}(u_a\pm v_a)$$ is an eigenvector of $\bH$ associated to  the eigenvalue $\pm s_a$.

 Then  we use the    classical trick by Erd\H{o}s, Schlein and Yau that for $E=\pm s_a$, $$|w_{a, \pm}(i)|^2\le \eta |G_\bH(E+\ii\eta)|_{ii}.$$ Then we prove that if $1\gg \eta\gg \max\{\eta_0,N^{-1/8}\},$ then  $$\|\E G_\bH(x+\ii\eta)\|=O(1),$$ uniformly on $E$ on $x\in [E\pm\eps]$, so that $$|\E G_\bH(E+\ii\eta)_{ii}|=O(1).$$ Then some concentration estimates allow to conclude.
 
   \section{Proof of Theorem \re{SRT} (Local Single Ring Theorem)}\la{secPrLSRT}
   
It is well known (see e.g. \cite[Sect. 4]{BOR1}) that   for any   $A\in \M$ and    any $F\in \Cc^2_c(\C)$,   
\be\la{251151}\int F(\lam)\ud \mu_A(\lam)=\ff{2\pi } \int_{z\in \C} \Del F(z)\lf(\int \log|s|\ud\nu_{z-A}(s)\ri)\ud \real(z)\ud \im(z).\ee

Here, we get 
\beq \int F_{z_0,\eps_N}(\lam)\ud \mu_A(\lam)&=& \ff{N}\sum_{i=1}^N F_{z_0,\eps_N}(\lam_i)\\
&=&\ff{2\pi}\int (\Del F_{z_0,\eps_N})(z)\int \log|s|\ud\nu_{z-A}(s)\ud \real(z)\ud \im(z)\\
&=&\ff{2\pi \eps_N^4}\int  (\Del f)(\f{z-z_0}{\eps_N})\int \log|s|\ud\nu_{z-A}(s)\ud \real(z)\ud \im(z)\\
&=&\ff{2\pi \eps_N^2}\int  \Del f(z')\int \log|s|\ud\nu_{z_0+\eps_Nz'-A}(s) \ud \real(z')\ud \im(z') 
\eeq
In the same way, we have (using an integration by parts at \eqre{2611513h} because $z_0$ has been chosen in the interior of   the support of the function $\rho$  of \eqre{bact15h20}),
\beqy\nonumber \int F_{z_0,\eps_N}(\lam)\ud\mu(\lam)&=& \ff{    \eps_N^2}\int f(\f{z-z_0}{\eps_N})\rho(z)\ud \real(z)\ud \im(z)\\ \la{2611513h}&=&\ff{2\pi  \eps_N^2}\int f(\f{z-z_0}{\eps_N})\Del_z(\int\log|x|\nu_{\infty,z}(\ud x))\ud \real(z)\ud \im(z)\\
\nonumber&=&\ff{2\pi  \eps_N^4}\int (\Delta f)(\f{z-z_0}{\eps_N}) \int\log|x|\nu_{\infty,z}(\ud x)\ud \real(z)\ud \im(z)\\
\nonumber&=&\ff{2\pi  \eps_N^2}\int (\Del f)(z') \int\log|x|\nu_{\infty,z_0+\eps_Nz'}(\ud x)  \ud \real(z')\ud \im(z')
\eeqy
Hence to prove that both expressions are equal up to an error term tending to zero in \pro as $N\to\infty$, by Lemma 3.1 of \cite{TaoVuKCL}, 
we need to prove that   
\bgt\ite[(i)]  for any $z'\in \C$, we have, for the convergence in probability,\bes\la{246141}
\eps_N^{-2}\lf|\int\log|x|\ud(\nu_{z_0+\eps_Nz'-A}-\nu_{\infty,z_0+\eps_Nz'}) (x)\ri| \Ninf 0,\ees 
\ite[(ii)]  for any $R>0$, the sequence \bes\la{246141bis}\eps_N^{-4}\int_{z'\in B(0,R)}\E\lf[ \lf|\int\log|x|\ud(\nu_{z_0+\eps_Nz'-A}-\nu_{\infty,z_0+\eps_Nz'}) (x)\ri|^2\ri]\ud\real(z')\ud\im(z')\ees  is bounded.
\ent
We shall in fact prove that for any $R>0$, uniformly in $z'\in B(0,R)$, $$\eps_N^{-4} \, \E\lf[ \lf|\int\log|x|\ud(\nu_{z_0+\eps_Nz'-A}-\nu_{\infty,z_0+\eps_Nz'}) (x)\ri|^2\ri] \Ninf 0,$$   which will prove (i) and (ii) in the same time. 



So let us fix $R>0$.

 Let us now choose a positive integer $p$ and $\epsilon>0$ \st \be\la{29111410h30minbf}4\al(p+2)+2\epsilon(p+1)<p\ee
 (which is possible since $4\al<1$) and set \be\la{29111410h0}t_N:=(\log N)^{-(2\al+\epsilon)}.\ee
Then \eqre{29111410h30minbf} implies that \be\la{29111410h}\f{t_N^{-(p+1)}}{(\log N)^{p/2}}\ll \eps_N^2\ee
and as $\epsilon>0$,    for any $k\ge 1$, \be\la{29111410h1}t_N|\log t_N|^k \ll \eps_N^2.\ee

Let 
 $\vfi_{t_N}$ be a smooth function with support contained in $[t_N/2, 3K+1]$, taking values in $[0,1]$,  equal to $1$ on $[t_N,3K]$. We can construct   a sequence of functions \st   there is a constant $C$ independent of $N$ \st for  all $\ell\ge 0 ,\ld, p+1$, \be\la{2611142}\|\vfi_{t_N}^{(\ell)}\|_\infty\le Ct_N^{-\ell}.\ee
We set \be\la{2611142SAMN}\log_{\ge t_N}(x):=\vfi_{t_N}(x)\log(x) \qquad\trm{ and }\qquad \log_{<t_N}(x):=(1-\vfi_{t_N}(x))\log(x).\ee
Then we have, for $\xi=z_0+\eps_Nz'$, 
\beqy\la{1-011} \lf|\int\log|x|\ud(\nu_{\xi-A}-\nu_{\infty,\xi}) (x)\ri|&\le&   \lf|\int\log_{\ge t_N}|x|\ud(\nu_{\xi-A}-\nu_{\infty,\xi}) ( x)\ri|\\ \nonumber&&\qquad +\lf|\int\log_{< t_N}|x|\ud \nu_{\xi-A} ( x)\ri|\\ \nonumber&&\qquad \qquad  +\lf|\int\log_{< t_N}|x|\ud \nu_{\infty,\xi} ( x)\ri| 
\eeqy
Let us treat the three terms in the RHS of \eqre{1-011} separately.

$\bullet$  
By Theorem \re{132151}, \eqre{2611142} and \eqre{29111410h},   with  probability at least $1-C\mre^{-N^c}$, uniformly in $z'\in B(0,R)$, we have 
$$\lf|\int\log_{\ge t_N}|x|\ud(\nu_{\xi-A}-\nu_{\infty,\xi}) ( x)\ri|
 \ \le \  C\f{t_N^{-(p+1)}}{(\log N)^{p/2}}.
$$ 
But by \eqre{29111410h}, the RHT is $\ll \eps_N^2$. As, on the complementary of the above event,  we   have the domination inequality $$\eps_N^{-4}\lf|\int\log_{\ge t_N}|x|\ud(\nu_{\xi-A}-\nu_{\infty,\xi}) ( x)\ri|^2\le 
4\log N(|\log t_N|+K)^2\ll C\mre^{N^c},$$ we deduce that unformly in $z'\in B(0,R)$,  \be\la{30111412h50}\E\lf[\eps_N^{-4}\lf|\int\log_{\ge t_N}|x|\ud(\nu_{\xi-A}-\nu_{\infty,\xi}) ( x)\ri|^2\ri]\Ninf 0.
\ee

$\bullet$ Let us now treat the close-to-zero terms. 
\beg{lem}\la{28111400283}Let $\mc{K}$ be a compact subset  of $\C$ which does not contain $0$.
\bgt\ite[(a)]Let $C$ be as in Hypothesis \re{2711141}, (iii). Then for any $\xi\in \C$,  $t\in (0,1)$, $$\lf|\int_{[0,t]} \log|x|(\nu_{\infty,\xi}) (\ud x)\ri|\le C t(1-\log(t)).$$  
\ite[(b)] There are some constants $A_{\mc{K}},a_{\mc{K}}>0$ \st for all $N$ large enough, for all $\xi\in \mc{K}$ and all $y>0$, $$\E[\nu_{\xi-A}([-y,y])] \ \le \ A_{\mc{K}}\max\{y, N^{-a_{\mc{K}}}\}.$$
\ite[(c)] 
Let  $N\ge 1$ and $\nu$ be a \pro measure on $\R$ \st for some constants $A,a>0$, for all $y>0$, \bes\la{27111414h}\nu ([-y, y])\le A\max\{y, N^{-a}\}.\ees Then there is $A'=A'(A,a)$   \st for any $t\in [N^{-a},1]$, we have  \bes\la{27111414h1}\int_{  [N^{-a},t]}|\log |x||^2\ud \nu (x)\le A'(t|\log t|^2+N^{-a}|\log N|^2).\ees
\ite[(d)] There are some positive   constants $c_{\mc{K}},C_{\mc{K}}$ \st for   $N\ge 1$ large enough, for  all $\xi\in \mc{K}$,  $u>0$, 
$$\p(s_{\min}(\xi-A)\le u)\le C_{\mc{K}}u^{c_{\mc{K}}}N^{C_{\mc{K}}}$$
and \st for any $\del>0$, we have, for all $\xi\in \mc{K}$,  $$ \E[|\log(s_{\min}(\xi-A))|^4\one_{s_{\min}(\xi-A)\le N^{-\del}}]\le C_{\mc{K}}N^{C_{\mc{K}}-\del c_{\mc{K}}}(\log N)^4.$$
\ent
\en{lem}

\bpr (a) For $C$ as in Hypothesis \re{2711141}, (iii), by Lemma \re{lemmeLipbxp}, we have   $\im m_\nu\le C$ on $\C^+$, hence   $\im m_{\nu^s}\le C$  on $\C^+$. So by Lemma \re{lemmeLipbxp}, for any $\xi\in \C$,  $\nu_{\infty, \xi}=\nu^s\bxp\del_{|\xi|}^s$ has a density with respect to the Lebesgue measure, which is bounded by $C/\pi$.  Thus for $t\in (0,1)$,  $$|\int_{[0,t]} \log|x|(\nu_{\infty, \xi}) (\ud x)|\le -(C/\pi)\int_0^{t}\log x\ud x=-(C/\pi)t(\log(t)-1).$$

(b) follows directly from Lemmas 13 and 15 of \cite{GUI} (the fact that the estimate is uniform in $\xi$ as $\xi$ stays bounded and bounded away from zero follows from a careful look at the arguments of \cite{GUI}).

(c) can be found in the proof of \cite[Prop. 4 (i), p. 1208]{GUI}. 

(d) The first part follows from Theorem 1.1 of \cite{RUD} by Rudelson and Vershynin, as for $\xi\ne 0$, $\ds s_{\min}(\xi-A)= |\xi|s_{\min}(U^*V^*-T/\xi)$.
Then to compute the   expectation, we will use the fact that 
for any positive random variable $X$, any $\al>0$ and any $\eps\in (0,1]$,  \bes\la{27111419h4}\E[|\log(X)|^\al\one_{X\le \eps}]=\al  \int_0^\eps\p(X\le   u)|\log u|^{\al-1}\f{\ud u}{u}+\p(X\le   \eps)|\log \eps|^\al,\ees
so that  
\beq \E[|\log(s_{\min}(\xi-A))|^4\one_{s_{\min}(\xi-A)\le N^{-\del}}]&\le &-2C_{\mc{K}}N^{C_{\mc{K}}} \int_0^{N^{-\del}} u^{c_{\mc{K}}-1}  \log^3 (u)  \ud u \\&& \qquad\qquad +C_{\mc{K}}c_{\mc{K}}\del N^{C_{\mc{K}}-c_{\mc{K}}\del}(\log N)^4\\
&\le &C'_{\mc{K}}N^{C_{\mc{K}}-\del c_{\mc{K}}}\lf(1 +(\log N)^4 \ri).
\eeq
\epr

We know that for $N$ large enough, for any $\xi\in B(z_0, R\eps_N)$, the support of $\nu_{\infty, \xi}=\nu^s\bxp\del_{|\xi|}^s$ and the spectrum of $|\xi-A|$ are contained in $[-3K,3K]$, so their intersection with the support of the function $\log_{<t_N} $ defined at \eqre{2611142SAMN}  is contained in $[0,t_N]$. As the function $\vfi_{t_N}$ only takes values in $[0,1]$, we have  \be\la{2811140029}\lf|\int\log_{< t_N}|x|\ud \nu_{\infty,\xi} ( x)\ri|\le -2\int_{[0,t_N]}\log x\ud \nu_{\infty,\xi} ( x) \ee and 
\beqy \nonumber\lf|\int\log_{< t_N}|x|\ud \nu_{\xi-A} ( x)\ri|& \le& -2\int_{[0,t_N]}\log x\ud \nu_{\xi-A} ( x) \\ \la{3011414h}&\le& -2\int_{[0,N^{-\del}]}\log x\ud \nu_{\xi-A} ( x)\\ \nonumber&&\qquad -2\int_{[N^{-\del}, N^{-a_{\mc{K}}}]}\log x\ud \nu_{\xi-A} ( x)\\ \nonumber&&\qquad \qquad -2\int_{[ N^{-a_{\mc{K}}}, t_N]}\log x\ud \nu_{\xi-A} ( x)\eeqy  
where $\del>0$ is chosen \st for $c_{\mc{K}}, C_{\mc{K}}$ is in (d) of the previous lemma, we have $C_{\mc{K}}-\del c_{\mc{K}}<0$.

$\bullet\bullet$ By \eqre{2811140029}, (a) of Lemma \re{28111400283} and \eqre{29111410h1}, we know that 
\be\la{281114002900}\lf|\int\log_{< t_N}|x|\ud \nu_{\infty,\xi} ( x)\ri|\ll \eps_N^2\ee

$\bullet\bullet$ Let us now treat the three terms of the RHS of \eqre{3011414h}.

\noindent{\bf First term  of the RHS of \eqre{3011414h}:} We always have  $$\int_{[0, N^{-\del}]}|\log |x||\ud \nu_{\xi-A} ( x)\le \nu_{\xi-A} ( [0, N^{-\del}]) |\log(s_{\min}(\xi-A))|\one_{s_{\min}(\xi-A)\le N^{-\del}}$$
Let us now take the second moment. 
By Cauchy-Schwartz, we have \beq \E\lf[\lf(\int_{[0, N^{-\del}]}|\log |x||\ud \nu_{\xi-A} ( x)\ri)^2\ri]
&\le & \E\lf[( \nu_{\xi-A} ( [0, N^{-\del}]))^2 |\log(s_{\min}(\xi-A))|^2\one_{s_{\min}(\xi-A)\le N^{-\del}}\ri]\\
&\le& \sqrt{\E[\nu_{\xi-A} ( [0, N^{-\del}])^4] \E[|\log(s_{\min}(\xi-A))|^4\one_{s_{\min}(\xi-A)\le N^{-c}}]}
\eeq
Then we use (b) of Lemma \re{28111400283} (plus the fact that $x^4\le x$ when $x\in [0,1]$) to upper bound $\E[\nu_{\xi-A} ( [0, N^{-\del}])^4]$ and (d) of Lemma \re{28111400283}  to upper bound $\E[|\log(s_{\min}(\xi-A))|^4\one_{s_{\min}(\xi-A)\le N^{-c}}]$. As we chose $\del$ so that $C_{\mc{K}}-\del c_{\mc{K}}<0$, we get that 
\bes\la{30111414h57}\E\lf[\lf(\int_{[0,N^{-\del}]}\log x\ud \nu_{\xi-A} ( x)\ri)^2\ri]\ll \eps_N^4.\ees

\noindent{\bf Second  term  of the RHS of \eqre{3011414h}:}  We have \bes\la{30111414h572}\E\lf[\lf|\int_{[N^{-\del}, N^{-a_{\mc{K}}}]}\log x\ud \nu_{\xi-A} ( x)\ri|^2\ri]\le \del\log N \E[\nu_{\xi-A} ( [N^{-\del}, N^{-a_{\mc{K}}}])]\le  C\del(\log N )N^{-a_{\mc{K}}}\ll \eps_N^4\ees where we used (b) of Lemma \re{28111400283}.

\noindent{\bf Third term  of the RHS of \eqre{3011414h}:}   By  (b) and (c) of Lemma \re{28111400283} and \eqre{29111410h1} (using Cauchy-Schwartz again, as above), we can claim that \bes\la{30111414h573}\E\lf[\lf(\int_{[ N^{-a_{\mc{K}}}, t_N]}\log x\ud \nu_{\xi-A} ( x)\ri)^2\ri]\ll \eps_N^4.\ees

$\bullet$ Let us conclude the proof.   By what precedes, we have proved that the RHS of \eqre{3011414h} has second moment $\ll \eps_N^4$, uniformly in $z'\in B(0,R)$. Besides, by 
 \eqre{281114002900}, we have proved that the (deterministic) RHT of \eqre{2811140029} is $\ll \eps_N^2$, uniformly in $z'\in B(0,R)$. This proves that the close-to-zero terms in \eqre{1-011}  have second moments $\ll\eps_N^{4}$,  uniformly in $z'\in B(0,R)$. At \eqre{30111412h50}, we proved that the same holds for the far-from-zero term in \eqre{1-011}. This concludes the proof of the theorem.

  \section{Proof of Theorem \re{241417h3NH} (Matrix Subordination)}\la{PrMaSub200115}

This proof goes roughly along the same lines as the one of Theorem 2 of the paper \cite{KarginAOP13Sub} by Kargin. The main difficulty is to give a Schwinger-Dyson equation adapted to our context (Lemma \re{164146h31lem}), which forces us to introduce the linear form $\tau$ of \eqre{174141} (from the point of view of quantum \pro theory, which identifies  the normalized trace to an expectation, $\tau$ can be assimilated to a conditional expectation).
 \subsection{Preliminaries}
    First, one can easily see, by left and right invariance of the Haar measure, that one can suppose that $A$ and $B$ are diagonal matrices with non negative entries, so that \be\la{941416h30}A^*=A\qquad ;\qquad B^*=B.\ee 
 
 To state our forthcoming Equation \eqre{164146h31},  we   define the map \beqy\la{174141} \tau : \mc{M}_{2N}(\C)= \bpm \M&\M\\ \M&\M\epm &\lto & 
   \bpm\M&0\\ 0&\M\epm\\  \nonumber\\
  \bpm A&B\\ C&D\epm &\longmapsto & 
  \bpm\ff{N}\Tr A&0\\ 0&\ff{N}\Tr D\epm,\nonumber\eeqy where the  complex numbers $\ff{N}\Tr A$ and $\ff{N}\Tr D$ are assimilated to the corresponding scalar matrices.

 \beg{rmk}\la{941416h301}  Let us introduce the matrix \be\la{defmatrixP416}P:=\bpm 0&I\\ I&0\epm\ee ($I$ denotes the identity matrix), which satisfies 
\be\la{defmatrixP4162}\forall X_1,X_2,Y_1,Y_2\in \mc{M}_N(\C), \qquad P\bpm X_1&Y_1\\ Y_2& X_2\epm P^{-1}= \bpm X_2&Y_2\\ Y_1& X_1\epm.\ee   Thus by   \eqre{941416h30},  $\bA$ and $\bB$ are invariant under conjugation by $P$ and 
  the matrix $\bW$ is   invariant, in law,  under    conjugation  by $P$. We deduce   that the (random or deterministic) matrices $\bH$, $G_\bH(z)$ and $\E G_\bH(z)$  are  invariant, in law,  under  conjugation by $P$. It implies  that $$\E[\tau(G_\bH(z))]\qquad;\qquad \E[\tau(G_\bH(z)\wt{\bB})]$$  are scalar $2N\ti 2N$ matrices equal to respectively $\E[m_\bH(z)]I$ 
 and $\E[f_B(z)]I$ for  \be\la{941416h11}f_B(z):=\ff{2N}\Tr (G_\bH(z)\wt{\bB}).\ee \en{rmk}

 The following lemma is the  \emph{Schwinger-Dyson equation} of  our problem.
 \beg{lem}\la{164146h31lem}For any $z$, 
  \be\la{164146h31}\E[\tau(G_{\bH}(z)) \wt{\bB}G_{\bH}(z)]=\E[\tau(G_{\bH}(z)\wt{\bB})  G_{\bH}(z)].\ee
 \en{lem}

 \bpr It suffices to prove that the element of $\Md\otimes \Md$  \be\la{841416h1}\E[G_{\bH}(z)\otimes (\wt{\bB}G_{\bH}(z))]-\E[(G_{\bH}(z)\wt{\bB})\otimes G_{\bH}(z)] \ee belongs to the kernel of the linear map $X\otimes Y\mapsto \tau(X)Y$. We shall  prove that  \eqre{841416h1} belongs to the space \be\la{841416h}\lf( \bpm \M&0\\ \M&0\epm\otimes \bpm 0&0\\ \M&\M\epm\ri)\oplus \lf(\bpm 0&\M\\ 0&\M\epm \otimes \bpm \M&\M\\ 0&0\epm\ri) \ee which is of course enough. Let us define 
 $$\Psi : \Md\otimes \Md\to \mc{L}(\Md)$$ to be the linear map defined by   $\Psi(X\otimes Y)(M)=XMY.$ It is easy to see that the space of 
  \eqre{841416h} is precisely the space of elements of $T\in \Md\otimes \Md$ \st  $$\bpm \M&0\\ 0&\M\epm\subset \ker \Psi(T).$$ Hence it suffices to prove that for any $Z,Z'\in \M$, \be\la{841416h2}\E[G_{\bH}(z)\bpm Z&0\\ 0& Z'\epm\wt{\bB}G_{\bH}(z)-G_{\bH}(z)\wt{\bB}\bpm Z&0\\ 0& Z'\epm G_{\bH}(z)]=0\ee By linearity, one can suppose that $Z,Z'$ are Hermitian. Then one recognizes easily that the LHT of \eqre{841416h2} is (up to a constant factor) the derivative, at $t=0$, of $\E[G_{\bH_t}(z)]$, where $$\bH_t:=\bA+\bpm \mre^{\ii t Z}&0\\ 0& \mre^{\ii tZ'}\epm \bW\bB\bW^*\bpm \mre^{-\ii t Z}&0\\ 0& \mre^{-\ii tZ'}\epm.$$ By invariance of the Haar measure, we have $$\bpm \mre^{\ii t Z}&0\\ 0& \mre^{\ii tZ'}\epm \bW\eqlaw \bW,$$ hence the above derivative is null. This proves the lemma.
 \epr
 
 Let   \be\la{941416h12}S_B(z):=-\f{\E[f_B(z)]}{\E[m_\bH(z)]}\ee for $f_B(z)$ and $m_\bH(z)$ defined at \eqre{941416h11} and \eqre{941419h4}.
 For $X,Y$ matrices, let $[X,Y]:=XY-YX$ denote the commutant of $X$ and $Y$.
 \beg{lem}\la{164146h41} Let $$\Del_A(z)=-\E[\oci{\tau(G_\bH(z))}G_\bH(z)]-\E[[\oci{\tau(G_\bH(z))},G_\bA(z)]\wt{\bB}G_\bH(z)]-\E[ \oci{\tau(G_\bH(z)\wt{\bB})} G_\bH(z)]$$ and \be\la{941416h12irmar}R_A(z):=\f{G_\bA(z+S_B(z))(\bA-z)\Del_A(z)}{\E[m_\bH(z)]}.\ee
 Then we have $$\E[G_\bH(z)]=G_\bA(z+S_\bB(z))+R_A(z).$$
 \en{lem}
 
 \bpr By Remark \re{941416h301}, $\E[\tau(G_\bH(z))]$ is a scalar $2N\ti 2N$ matrix  equal to $\E[m_\bH(z)]I$.
 Thus, using successively  \eqre{941416h33}, the resolvant identity and the previous lemma,  we get \beq\E[m_\bH(z)]\E[G_\bH(z)]&=&\E[\tau(G_\bH(z))]\E[G_\bH(z)]\\
 &=&\E[\tau(G_\bH(z))G_\bH(z)]-\underbrace{\E[\oci{\tau(G_\bH(z))}G_\bH(z)]}_{\ds :=\eps_1}\\
 &=&\E[\tau(G_\bH(z))\{G_\bA(z)-G_\bA(z)\wt{\bB}G_\bH(z)\}]-\eps_1\\
  &=&\E[m_\bH(z)]G_\bA(z)-\E[\tau(G_\bH(z))G_\bA(z)\wt{\bB}G_\bH(z)]-\eps_1\\
  &=&\E[m_\bH(z)]G_\bA(z)-\underbrace{\E[[\tau(G_\bH(z)),G_\bA(z)]\wt{\bB}G_\bH(z)]}_{\ds :=\eps_2}\\ &&-\E[G_\bA(z)\tau(G_\bH(z))\wt{\bB}G_\bH(z)]-\eps_1\\
    &=&\E[m_\bH(z)]G_\bA(z)-G_\bA(z)\E[ \tau(G_\bH(z)\wt{\bB})G_\bH(z)]-\eps_1-\eps_2\\
     &=&\E[m_\bH(z)]G_\bA(z)-G_\bA(z)\E[ \tau(G_\bH(z)\wt{\bB})]\E[G_\bH(z)]\\
     &&-\underbrace{G_\bA(z)\E[ \oci{\tau(G_\bH(z)\wt{\bB})} G_\bH(z)]}_{\ds :=\eps_3}-\eps_1-\eps_2\\
     &=&\E[m_\bH(z)]G_\bA(z)-G_\bA(z)\E[ \tau(G_\bH(z)\wt{\bB})]\E[G_\bH(z)]-\eps_1-\eps_2-\eps_3
     \\
     &=&\E[m_\bH(z)]G_\bA(z)-G_\bA(z)\E[ f_B(z)]\E[G_\bH(z)]-\eps_1-\eps_2-\eps_3
 \eeq
 Dividing by the complex number $\E[m_\bH(z)]$ and multiplying on the left by $\bA-z$, one gets $$(\bA-z)\E[G\bH(z)]=I+S_B(z)\E[G_\bH(z)]-\eps'_1-\eps'_2-\eps'_3,$$ for $\ds \eps_i':=\f{(\bA-z)\eps_i}{\E[m_\bH(z)]}$.
This gives $$(\bA-z-S_B(z))\E[G_\bH(z)]=I-\eps'_1-\eps'_2-\eps'_3,$$ \ie $$ \E[G_\bH(z)]=G_\bA(z+S_B(z))-\eps''_1-\eps''_2-\eps''_3,$$ for $$\ds \eps_i'':=G_\bA(z+S_B(z))\eps_i'=\f{G_\bA(z+S_B(z))(\bA-z)\eps_i}{\E[m_\bH(z)]}.$$
To conclude, it suffices to notice that $$R_A(z)=-\eps''_1-\eps''_2-\eps''_3,$$
 up to the fact that in the second term of $R_A(z)$, we have $[\oci{\tau(G_\bH(z))},G_\bA(z)]$ instead of $[\tau(G_\bH(z)),G_\bA(z)]$. But as $\E[\tau(G_\bH(z))]$ is a scalar matrix, both are equal.
 \epr
 
 \beg{lem}\la{441410h}Let $\Psi_A(z):=\ff{\E m_\bH(z)}(\bA-z) \E[\Del_A(z)]$ and $Y_A(z):=(I+\Psi_A(z))^{-1}-I$. Then \be\la{441416h1}S_B(z)I=-(\E G_\bH(z))^{-1}+\bA-z+Y_A(z)(\bA-z-S_B(z)),\ee where $I$ denotes the identity matrix.
 \en{lem}
 
 \bpr By the previous lemma,   \beq \E G_\bH(z)&=& G_\bA(z+S_B(z))\lf(I+\ff{\E m_\bH(z)}(\bA-z) \E\Del_A\ri)\\ &=& G_\bA(z+S_B(z))\lf(I+\Psi_A(z)\ri)\\
 &=&G_\bA(z+S_B(z))\lf(I+Y_A(z)\ri)^{-1}
 \eeq
  hence $$(\E G_\bH(z))^{-1}=\lf(I+Y_A(z)\ri)(\bA-(z+S_B(z)))$$
 which allows to conclude.
 \epr
 
  \beg{lem}Let $\rho$ be a \pro measure supported by $[-K,K]$ and $|z|\le K$. Then  $$\im m_\rho(z) \ge \f{\eta}{5K^2}.$$  \en{lem}
 
 \bpr It suffices to note that for any   $\lam\in [-K,K]$, $\ds\im \ff{\lam-z}=\f{\eta}{(\lam-E)^2+\eta^2}\ge \f{\eta}{(2K)^2+K^2}.$
 \epr
 
 It follows from this lemma that  there is $c$ depending only on $K$ \st  \be\la{41416h2MLB}\ff{|\E m_\bH(z)|}\le \f{c}{\eta}\ee and \be\la{41416h2}|S_B(z)|=\lf|\f{\E f_B(z)}{\E m_\bH(z)}\ri|\le \f{c}{\eta^2}.\ee

  \beg{lem}For any $\del>0$, $$\p(\|\oci{\tau(G_\bH(z))}G_\bH(z)\|\ge \del ) \ \le \ 4\exp\lf(-c\f{\del^2\eta^6}{\|B\|^2}N^2\ri)$$
 $$\p(\|[\oci{\tau(G_\bH(z))}, G_\bA(z)]\wt{\bB}G_\bH(z)\|\ge \del ) \ \le \ 4\exp\lf(-c\f{\del^2\eta^8}{\|B\|^4}N^2\ri)$$ and $$\p(\|\oci{\tau(G_\bH(z)\wt{\bB})}G_\bH(z)\|\ge \del )\ \le\ 4 \exp\lf(-c\f{\del^2\eta^6}{\|B\|^4(1+\f{\eta}{\|B\|})^2}N^2\ri)$$
 \en{lem}
 
\bpr To prove it, as  for any Hermitian matrix $M$, $\|G_M(z)\|\le \eta^{-1}$, it suffices to prove that 
 for any $\del>0$, $$\p(\|\oci{\tau(G_\bH(z))}\|\ge \del) \ \le \ 4\exp\lf(-c\f{\del^2\eta^4}{\|B\|^2}N^2\ri)$$ and $$\p(\|\oci{\tau(G_\bH(z)\wt{\bB})}\|\ge \del)\ \le\ 4 \exp\lf(-c\f{\del^2\eta^4}{\|B\|^4(1+\f{\eta}{\|B\|})^2}N^2\ri).$$

 We shall apply the   Lemma \re{1811141} of the appendix. Note first that   a Haar-distributed unitary matrix can be realized as the product of a Haar-distributed $\mc{SU}_N$ matrix by a uniform random phase, hence up to a randomization of  $B$ by multiplication by an independent  uniform phase,   one can suppose that $U$ and $V$ are independent  Haar-distributed $\mc{SU}_N$ matrices.

Let $P_1, P_2$ be the $2N\ti 2N$ matrices defined by $$P_1:=\bpm I&0\\ 0&0\epm\qquad 
P_2:=\bpm 0&0\\ 0&I\epm,$$ so that for any $M\in \Md$, $$\tau(M)=\bpm N^{-1}\Tr P_1MP_1&0\\ 0&N^{-1}\Tr P_2MP_2\epm.$$

Let $\vfi_i$,   $\psi_i$ ($i=1,2$)    be the functions defined on $(\mc{SU}_N)^2$ by $$\vfi_i(U,V):=N^{-1}\Tr (P_iG_\bH(z)P_i)\qquad \psi_i(U,V):=
N^{-1}\Tr (P_iG_\bH(z)\wt{\bB}P_i)$$ with the notations of \eqre{1141413h}, \eqre{1141413h1}. We need to prove that under the sole hypothesis that $\|A\|, \|B\|, |z|\le K$, the numbers $$\f{NL_{\vfi_i}^2\eta^4}{\|B\|^2}\qquad;\qquad \f{NL_{\psi_i}^2\eta^4}{\|B\|^4 (1+\f{\eta}{\| B\|})^2}\qquad ; \qquad (i=1,2)
$$ are bounded uniformly in $N$. 

For $X,Y$  skew-Hermitian matrices with null traces and $U,V \in \mc{SU}_N$, $$\pa_{t,t=0}\vfi_i(\mre^{tX}U,\mre^{tY}V)=-N^{-1}\Tr (P_iG_\bH(z)(\bZ\wt{\bB}-\wt{\bB}\bZ)G_{\bH}(z)P_i),$$ with $\bZ:=\bpm X&0\\ 0&Y\epm $, so that  $$\nabla \vfi_i(U,V)=-\ff{N} \bW^*\mc{P}( [\wt{\bB},G_\bH(z)P_iG_\bH(z)]),$$ where $\mc{P}$ is the orthogonal projection from  $\Md$ onto the tangent space at $I$, of $(\mc{SU}_N)^2$. As this projection does not enlarge the norm, the usual non commutative H\H{o}lder inequalities (see Appendix A.3 of \cite{agz}) allow to claim that  $$\f{NL_{\vfi_i}^2\eta^4}{\|B\|^2}$$ is bounded.

In the same way, 
$$\pa_{t,t=0}\psi_i(\mre^{tX}U,\mre^{tY}V)=N^{-1}\Tr \{P_i(-G_\bH(z)(\bZ\wt{\bB}-\wt{\bB}\bZ)G_{\bH}(z)\wt{\bB}+ G_\bH(z)(\bZ\wt{\bB}-\wt{\bB}\bZ))P_i\},$$
hence  $$\nabla \vfi_i(U,V)=\ff{N}W^*\mc{P}([\wt{\bB},-G_\bH(z)\wt{\bB}P_iG_\bH(z)+P_iG_\bH(z)]),$$ and one concludes as above.  
\epr

  \beg{lem}\la{lem171142148}We have $$\E\|\oci{\tau(G_\bH(z))}G_\bH(z)\|  \ \le \ c\f{\|B\|}{N\eta^3}$$
 $$\E\|[\oci{\tau(G_\bH(z))}, G_\bA(z)]\wt{\bB}G_\bH(z)\|  \ \le \ c\f{\|B\|^2}{N\eta^4}$$ and $$ \E\|\oci{\tau(G_\bH(z)\wt{\bB})}G_\bH(z)\| \ \le\ c\f{\|B\|^2(1+\f{\eta}{\|B\|})}{N\eta^3}$$
 Hence if $|z|, \|B\|\le K$, then \be\la{241417h2BOB}\|\Del_A(z)\| \le c\ff{N\eta^4}\qquad;\qquad  \|\Psi_A(z)\|\le  c\ff{N\eta^5}\ee and  \be\la{441416h3}N \eta^{5}\ge 2c\implies \|Y_A(z)\|\le \f{2c}{N\eta^5}.\ee
 \en{lem}
 
 \bpr
 The three first inequalities follow from the previous lemma and standard queues-moments relations. The upper-bound on $\|\Del_A(z)\| $ follows from the very definition of $\Del_A(z)$ at Lemma \re{164146h41}.
The upper-bound on $\|\Psi_A(z)\| $  follows from its definition $$\Psi_A(z)=\ff{\E m_\bH(z)}(\bA-z) \E[\Del_A(z)]$$ and from  \eqre{41416h2MLB}. At last, \eqre{441416h3} follows from the definition $Y_A(z)=(I+\Psi_A(z))^{-1}-I$ and from the well known inequality  $$\|X\|\le   1/2\implies \|(I-X)^{-1}-I\|\le 2\|X\|.$$  \epr
 
 Adapting the proof  of Lemma 4.7 of \cite{bbcf12}, we get the following lemma.

\beg{lem}Let $\U\subset \mc{M}_N(\C)$ be a compact  Lie group  and let $\mc{L}$ be the complex linear subspace of $\mc{M}_N(\C)$ spanned by its Lie algebra.  Let us fix $M\in \mc{M}_N(\C)$ and, for $b\in \mc{M}_N(\C)$ \st $b-uMu^{-1}$ is invertible  for any $u\in \U$ and define the random matrix  $R(b):=(b-UMU^{-1})^{-1}$, where $U$ is Haar-distributed in $\U$. 
Then:  \bgt\ite[(i)]
for any $Y\in \mc{L}$, we have $\E[R(b)]Y-Y\E[R(b)]=\E[R(b)(Yb-bY)R(b)]$,
\ite[(ii)] the matrix $\E[R(b)]$ commutes with any matrix in $\mc{L}$  commuting with $b$. 
\ent \en{lem}

\bpr (ii) is a direct consequence of (i). Besides, by linearity, it suffices to prove (i) for $Y$ in the Lie algebra of $\U$. For such a matrix $Y$, differentiating at $0$ the constant function $f(t):=\E[(b-\mre^{tY}UMU^{-1}\mre^{-tY})^{-1}]$, we get $$\E[R(b)(UMU^{-1}Y- YUMU^{-1})R(b)]=0$$
Then, using that $R(b)UMU^{-1}=-I+R(b)b$ and that $UBU^{-1}R(b)=-I+bR(b)$, we get (i) directly.  
\epr

 \beg{lem}\la{164146h45} The matrix $\E G_\bH(z) $   commutes with $\bA$.\en{lem}
 
 \bpr Let us apply the previous lemma for $\U$ the group   of matrices  
$\bpm U& 0\\ 0& V\epm$, for $U,V$   unitary matrices (so that  $\mc{L}$ is the space of matrices   $\bpm X& 0\\ 0& X'\epm$, for $X,X'\in \mc{M}_N(\C)$), $M=\bB  $   and $b=z-\bA$. It states  that  $\E G_\bH(z) $   commutes with any matrix of $\mc{L}$ commuting with $\bA$, for example   with  any matrix of the type $\bpm D& 0\\ 0& D\epm$, with $D$ diagonal. We deduce that  $$\E G_\bH(z) =\bpm J&K\\ L&M\epm,$$ with $J,K,L,M$ some $N\ti N$ diagonal matrices. But by Remark \re{941416h301}, the matrix $\E G_\bH(z) $ is invariant under conjugation by the matrix $P$ introduced at \eqre{defmatrixP416}. By \eqre{defmatrixP4162}, it implies  that $J=M$ and $K=L$. It suffices  to conclude.
 \epr

 \subsection{Proof of Theorem \re{241417h3NH} (Matrix Subordination)}
 Note first that the statement  is   symmetric in $A$ and $B$, so we shall prove it for $S_B$ and $R_A$ only.  
 
 By Lemmas \re{164146h41} and \re{441410h}, 
   we have $$\E[G_\bH(z)]=G_\bA(z+S_\bB(z))+R_A(z)$$
with  \be\la{241417h2BB}R_A(z):=\f{G_\bA(z+S_B(z))(\bA-z)\Del_A(z)}{\E[m_\bH(z)]}\ee
and $$S_B(z)I=-(\E G_\bH(z))^{-1}+\bA-z+Y_A(z)(\bA-z-S_B(z)).$$
 By Proposition 2 of \cite{KarginAOP13}, we know that $$-(\E G_\bH(z))^{-1}$$ has all its eigenvalues with imaginary part $\ge \eta$. But   by Lemma \re{164146h45}, $\E G_\bH(z)$ commutes with $\bA$, hence  the eigenvalues of $-(\E G_\bH(z))^{-1}+\bA-z$ have non  negative imaginary parts.  Thus by standard perturbation analysis (see e.g. \cite[Chap. 4]{CGLP}), $$\im S_B(z)\ge  \|Y_A(z)(\bA-z-S_B(z))\|.$$
 Then, \eqre{41416h2} and  \eqre{441416h3} give directly the lower-bound \eqre{241417h} on the imaginary part of $S_B(z)$.

 The upper-bound \eqre{241417h2} on $\|R_A(z)\|$ follows directly from the expression \eqre{241417h2BB}, the upper-bound \eqre{241417h2BOB} on $\|\Del_A(z)\|$, and the fact that  for $N$ large enough, $\im (z+S_B(z))\ge \f{\eta}{2}$.

 \section{Appendix}
 
 \subsection{Free convolution and subordination}
  Let us first recall one of the ways to define the free convolution \cite{BianePFE,BA,BB07}.   \beg{Th}[Definition of the free convolution via subordination]\la{24141}Let $\mu, \nu$ be \pro measures on the real line with compact supports. Then the system \beq\la{74141} m(z)&=& m_{\mu}(z+S_\nu(z))\\
 \la{74142}m(z)&=& m_{\nu}(z+S_\mu(z))\\ 
\la{74143} -\lf(z+\ff{m(z)}\ri)&=& S_\mu(z)+S_\nu(z)\eeq
 has a unique solution $(m(\cdot), S_\mu(\cdot), S_\nu(\cdot))$
   in the class of triplets of analytic functions on $\C^+$ satisfying, as $|z|\to+\infty$,\beqy\nonumber m(z)&=&-z^{-1}+O(z^{-2})\\\nonumber \\ \la{2311512h45}|S_\mu(z)|+|S_\nu(z)|&=&O(1)\eeqy
   The function $m(z)$ is then the Stieltjes transform of a unique \pro measure, which is $\mu\bxp\nu$. Moreover, $S_\mu$ and $S_\nu$ take values in $\ovl{\C^+}$. 
 \en{Th}
 \beg{rmk}\la{201152}Note that this result, in addition to define the free convolution, is a first \emph{regularity result} for this convolution. Indeed, for any $\eta>0$ and any \pro measure $\rho$ on the real line, the function $z\mapsto m_\rho(z+\ii\eta)$ is the Stieltjes transform of an analytic regularization of $\rho$ (namely its classical convolution with the Cauchy law $\ff{\pi}\f{\eta\ud x}{x^2+\eta^2}$). Hence for $\mu,\rho$ some \pro measures on the real line, the equation $$m_\rho(z)=m_\mu(z')\quad\trm{ with }\im z'> \im z$$ implies roughly that $\rho$ is more regular than $\mu$.
 \en{rmk}

The following lemmas will be used in this text.
 \beg{lem}\la{1871414:59}Let $X, Y$ be free self-adjoint elements of a tracial $W^*$-\pro space $ (\mc{A}, \tau)$ with repsective distributions $\mu, \nu$. Then for any $z\in \C^+$, we have $$S_\mu(z)=-\f{\tau(X(X+Y-z)^{-1})}{\tau((X+Y-z)^{-1})} \qquad;\qquad S_\nu(z)=-\f{\tau(Y(X+Y-z)^{-1})}{\tau((X+Y-z)^{-1})}.$$
\en{lem}

\bpr Let us focus for example on $S_\mu$. It is equivalent to prove that $$z+S_\mu(z)=z-\f{\tau(X(X+Y-z)^{-1})}{\tau((X+Y-z)^{-1})} =\f{\tau((z-X)(X+Y-z)^{-1})}{\tau((X+Y-z)^{-1})} $$ \ie that \be\la{1871413h55}(z+S_\mu(z))\tau((X+Y-z)^{-1})=\tau((z-X)(X+Y-z)^{-1}).\ee Let $\tau_{Y}$ denote the conditional (non-commutative) expectation given the $W^*$-algebra generated by $Y$.
We know, by Th. 3.1 of \cite{BianePFE}, that $$\tau_{Y}(\ff{z-X-Y})=\ff{z+S_\mu(z)-Y},$$ so that $$z+S_\mu(z)=Y+(\tau_{Y}(\ff{z-X-Y}))^{-1}$$ (the miracle of  \cite{BianePFE} being precisely that despite the $\tau_{Y}$ in  the RHT,  $z+S_\mu(z)$ is a scalar) and $$(z+S_\mu(z))\tau_{Y}(\ff{z-X-Y})=Y\tau_{Y}(\ff{z-X-Y})+1.$$ Let us now apply $\tau$. As $z+S_\mu(z)\in \C$,  we get 
$$(z+S_\mu(z))\tau((z-X-Y)^{-1})=\tau(Y(z-X-Y)^{-1})+1, $$ which is exactly \eqre{1871413h55}.
\epr

\beg{lem}\la{231151}
Let $K>0$ be fixed. Then there are $M_i=M_i(K)>0$ ($i=1,2$) \st for any pair $\mu, \nu$ of \pro measures with supports contained in $[-K,K]$, for any $|z|>M_1$, the numbers $S_{\mu^s}(z)$, $S_{\nu^s}(z)$ (defined at Theorem \re{24141}), $\ka_{\mu^s,\nu^s}(z)$,  $\al_{\mu^s,\nu^s}(z)$ and $\bet_{\mu^s,\nu^s}(z)$ (defined at \eqre{14111423h}, \eqre{defalpha} and \eqre{defbeta}) satisfy $$\ka_{\mu^s,\nu^s}(z)\ge \ff{M_2}\qquad ;\qquad \ff{\al_{\mu^s,\nu^s}(z) \bet_{\mu^s,\nu^s}(z)}\le M_2 \qquad ;\qquad \f{\al_{\mu^s,\nu^s}(z)}{|z+S_{\mu^s}(z)+S_{\nu^s}(z)|^2}\le M_2.$$
\end{lem}

\bpr First, by the previous lemma, we know that the estimate \eqre{2311512h45} is uniform in all pairs $\mu, \nu$ of \pro measures with supports contained in $[-K,K]$. Besides, it is obvious, from the series expansion,  that the estimates $-z^km_\mu^{(k)}(z)\lto 1$, $k=0,1,2$,  as $|z|\lto \infty$, are uniform in \pro measures $\mu$   with support contained in $[-K,K]$. Then, going back to the formulas defining  the functions of interest here, we get the desired estimates.
\epr

 \subsection{Concentration of measure for the Haar measure}

  By the lemma of  Gromov and Milman (see for example    \cite{agz}, page 299) and Proposition 1.11 of \cite{ledoux-amsbook}, we have:
\beg{lem}\la{1811141}Let $f$ be a smooth real-valued function on $(\mc{SU}_N)^2$ and let $$L_f:=\max_{(U,V)\in (\mc{SU}_N)^2}\sqrt{\Tr(\nabla f(U,V)\nabla f(U,V)^*)}.$$Then  for $U,V$     independent  Haar-distributed $\mc{SU}_N$ matrices, for any $\del\ge 0$, $$\p(|f(U,V)-\E[f(U,V)]|\ge \del)\le  2\exp\lf(-\f{N\del^2}{4L_f^2}\ri).$$ \en{lem}

 \subsection{Kantorovich's Theorem on Newton's Method}
 
 Let us give the simplified version of Kantorovich's Theorem that we need (particular case of \cite[Th. 1]{KantoTheo}). We let $\|\cdot\|$ denote the canonical euclidian norm on $\R^d$ or the associated  operator norm on $\mc{L}(\R^d)$.
 
 \beg{Th}\la{kanto}Let $O$ be an   open subset of $\R^d$ and $F: O\to \R^d$ be a $\Cc^1$ function. Let $x_0\in O$ \st $F'(x_0)$ is invertible and $y_0\in \R^d$. Suppose that for   $$L:=\sup_{x\ne y \in O}\f{\|F'(x_0)^{-1} (F'(x)-F'(y))\|}{\| x-y\|}\qquad;\qquad b:=\|F'(x_0)^{-1}(F(x_0)-y_0)\|$$
 we have $2bL<1.$
 Define $$r_*:= \f{2b}{1+\sqrt{1-2bL}}\le 2b\qquad ;\qquad r_{**}:=\f{1+\sqrt{1-2bL}}{L}$$ and choose $\rho\in [r_*, r_{**})$ \st $\ovl{B}(x_0, \rho)\subset O$. 
 Then the equation $F(x)=y_0$ has a unique solution $x_*$ in  $B(x_0, \rho) $ and this solution satisfies $$\|x_*-x_0\|\le r_*.$$
 \en{Th}
 
 
 \subsection{Local laws and Stieltjes transforms}

 \subsubsection{Density and upper bound on the Stieltjes transform}
  \beg{lem}\la{lemmeLipbxp}For $\mu$ \pro measure on the real line and $M>0$, we have equivalence between:\bgt\ite[(i)] $\mu$ admits a density $\rho$ with respect to the Lebesgue measure \st $\|\rho\|_\infty\le M$,
 \ite[(ii)] $\im m_\mu$ is uniformly bounded, on $\C^+$, by $\pi M$.
 \ent
 Moreover, in this case, for any \pro measure $\nu$ on $\R$, $\mu\bxp\nu$ also admits a density  with respect to the Lebesgue measure  which is bounded by  $M$.
 \en{lem}
 
 \bpr   For $z=E+\ii\eta$ ($E\in \R, \eta>0$), we have  $$ \ff{\pi}\im m_\mu(z)=\int_{\lam\in \R}\f{\eta}{(E-\lam)^2+\eta^2}\ud\mu(\lam).$$
 Hence if (i) holds, then $$  \im m_\mu(z)=\int_{\lam\in \R}\f{\eta}{(E-\lam)^2+\eta^2}\rho(\lam)\ud\lam\le M\int_{\lam\in \R}\f{\eta}{(E-\lam)^2+\eta^2} \ud\lam\le M\pi.$$
 Reciprocally, let us suppose that $\im m_\mu$ is uniformly bounded, on $\C^+$, by $\pi M$. The law $$\ff{\pi}\im m_\mu(\lam+\ii\eta)\ud \lam=\mu*\f{\eta\ud x}{\pi(x^2+\eta^2)} $$ converges weakly to $\mu$ as $\eta \da 0$, hence for any $x<y$, $$\mu([x,y])=\lim_{\eta\to 0}\ff{\pi}\int_x^y\im m_\mu(\lam+\ii\eta)\ud \lam$$ (indeed, by \cite[Lem. 2.17, (2)]{BA}, $\mu$ has no atom), so for any $x<y$, $\mu([x,y])\le M(y-x).$ This implies that the cumulative distribution function of $\mu$ is $M$-Lipschitz, hence is almost everywhere differentiable, with derivative $\le M$ and is the integral of its derivative, which is exactly (i).  
 
 The last statement follows from the   subordination for $\bxp$: by Theorem \re{24141},  there is a function $S:\C^+\to\ovl{\C^+}$ \st  on $\C^+$, $m_{\mu\bxp\nu}(z)=m_\mu(z+S(z))$, which allows directly to conclude by what precedes.
 \epr

 \subsubsection{Erd\H{o}s, Schlein and Yau's method}
   \beg{Th}[Erd\H{o}s, Schlein, Yau]\la{methodErdosSchleinYau2009}Let $\nu$ be a signed  measure on $\R$, $E\in \R$ and $\eta, M>0$.
   Then $\ds \f{\lf|\nu([E\pm M\eta])\ri|}{2M\eta}$ is upper-bounded by  \begin{align*}&C\bigg( \| m_{\nu}(\,\cdot\,+\mathrm{i}\eta) \|_{\infty,[E\pm M\eta]} +\f{|\nu|([E\pm 2M\eta])}{M^{3/2}\eta}+\\ &\qquad\qquad \f{|\nu|([E- 2M\eta\pm \sqrt{M}\eta]\cup [E+ 2M\eta\pm \sqrt{M}\eta])}{M\eta}+\f{\im m_{|\nu|}(E+\mathrm{i}M\eta) }{M}\bigg) ,\end{align*}
   for a certain universal constant $C$. 
   \en{Th}
   
   \bpr Let us briefly present   the ideas of the proof of Corollary 4.2 in \cite{ESY2}. It suffices to notice that for $\ds R(\lam):=\ff{\pi}\int_{[E\pm M\eta]}\f{\eta}{(x-\lam)^2+\eta^2}\ud x$, there are some functions  $T_1,T_2,T_3$ \st  
   $R(\lam)=\one_{|\lam-E|\le \eta M}+T_1(\lam)+T_2(\lam)+T_3(\lam),$ with \bgt \ite $\ds \|T_1\|_\infty\le \f{c}{\sqrt{M}}$ and $\supp(T_1)\subset [E\pm 2M\eta],$\ite 
   $\ds \|T_2\|_\infty\le 1$ and $\supp(T_2)\subset [E- 2M\eta\pm \sqrt{M}\eta]\cup [E+ 2M\eta\pm \sqrt{M}\eta],$
 \ite $\ds |T_3(\lam)|\le \f{CM\eta^2}{(\lam-E)^2+M^2\eta^2} .$\ent
 Hence as   $\ds\int_{\lam\in \R} R (\lam)\ud\nu(\lam)= \ff{\pi}\int_{[E\pm M\eta]}\im m_\nu(x+\mathrm{i}\eta)\ud x$, we have
\beq \nu ([E\pm M\eta]) 
   &=&\ff{\pi}\int_{[E\pm M\eta]}\im m_\nu(x+i\eta) \ud x
   -\int_{\lam\in \R}  (T_1(\lam)+T_2(\lam)+T_3(\lam))\ud \nu(\lam),\eeq which proves the theorem.
    \epr 
    

    \subsubsection{Helffer-Sj\H{o}strand functional calculus} The use of this method in random matrix theory   is quite recent  (see \cite[Proof of Lem. 5.5.5]{agz} or \cite[Proof of Lem. B.1]{ESRY}). As we shall use it in a non common scale (see Corollary \re{2611141}), we state it precisely here.

\beg{Th}\la{Helffer-Sjostrand}Let $\nu$ be a   signed measure on $\R$  and   $\phi$   a $\Cc^{p+1}$ compactly supported function on the real line, for $p\ge 1$.  Suppose that there are $\eta_{\min},\del>0$, $\al\in [0, p+1) $ \st   \be\la{18111417h}(E\in \op{supp}(\phi)\trm{ and }\eta>\eta_{\min})\implies    |m_{\nu}(E+\mathrm{i}\eta)|\le \del\eta^{-\al}. \ee
    Then $$\lf|\int\phi(x)\nu(\ud x)\ri|\le \f{L\|\phi^{(p+1)}\|_\infty}{p!\pi}|\nu|(\R) \lf(\f{\eta_{\min}^p}{p}+\f{\del  }{p-\al+1}\ri)  ,$$ for $L$ the Lebesgue measure of the support of $\phi$.
     \en{Th}
     
     \bpr 
     $\bullet$ The function \be\la{6121418h}f_{\nu}:(E, \eta)\in \R\ti[0,+\infty)\longmapsto  \beg{cases}\eta \, m_\nu(E+\mathrm{i}\eta) &\trm{ if $\eta>0$,}\\ 
\mathrm{i} \nu(\{E\})&\trm{ if $\eta=0$,}\\ \en{cases}
\ee   satisfies $\| f_\nu\|_{\infty}\le |\nu|(\R)$.
Indeed, for $\eta>0$, we have $$|f_\nu(E,\eta)|=\lf| \int \f{\ud \nu(x)}{x-\f{E}{\eta}+\ii\eta}\ri|\le \int {\ud |\nu|(x)} .$$
     
     $\bullet$ Choose   $\vfi:\R\to[0,1]$ a smooth function   with value $1$ in a neighborhood of  $0$ and supported by $[-1,1]$. Then set  $\ds\Psi(x+\mathrm{i}y):=\sum_{\ell=0}^p\f{\mathrm{i}^\ell}{\ell !}\phi^{(\ell)}(x)\vfi(y)y^\ell.$ Note that for $\ovl{\pa}:=  \pa_x+\mathrm{i}\pa_y $,  the functions  \be\la{6121418h01}  (E, \eta)\in \R\ti[0,+\infty)\longmapsto  \beg{cases}\ds \f{\ovl{\pa} \Psi(E+\mathrm{i}\eta)}{\eta^p}&\trm{ if $\eta>0$,}\\ 
 \f{\mathrm{i}^p}{p! }\phi^{(p+1)}(E)\vfi(0)&\trm{ if $\eta=0$,}\\ \en{cases}\ee      
  is     continuous and bounded.
Indeed, we have 
\be\la{6121418h25ccc}  \ovl{\pa} \Psi(E+\mathrm{i}\eta)\; =\;\f{\mathrm{i}^p}{p!}\phi^{(p+1)}(E)\vfi(\eta)\eta^p\ee  so that  the continuity is obvious and  \beqy  \sup_{E\in \R, \eta>0}\lf|\f{\ovl{\pa} \Psi(E+\mathrm{i}\eta)}{\eta^p}\ri|  &\le &\ff{p!}  \|\phi^{(p+1)}\|_\infty 
  \la{6121418h25}.\eeqy

$\bullet$ Using this remark in the particular case where $p=1$, we get     
\be\la{1061412h}\int_{t\in \R}\Psi(t)\ud\nu(t)=\pi^{-1}\real\lf(\iint_{(E, \eta)\in \R\ti[0,+\infty)}   \f{\ovl{\pa}\Psi(E+\mathrm{i}\eta)}{\eta}\eta m_\nu(E+\mathrm{i}\eta)\ud E\ud \eta\ri).\ee
Indeed,  \eqre{1061412h} is continuous (for the topology defined by bounded continuous   functions) and  linear in $\nu$, so that it suffices to prove it for  $\nu=\del_\lam$, with $\lam\in \R$. Then it is the content of \cite[Prop. C.1]{FBGAK2016}. 
  
  $\bullet$ As a consequence, as $\phi$ and $\Psi$ coincide on $\R$, using   \eqre{1061412h}, \eqre{18111417h} and \eqre{6121418h25},  we get    \beq \lf|\int\phi(x)\nu(\ud x)\ri|&\le&\pi^{-1}\lf|\iint_{(E, \eta)\in \supp(\phi)\ti (0,1]}   \f{\ovl{\pa}\Psi(E+\mathrm{i}\eta)}{\eta^p}\eta^p m_{\nu}(E+\mathrm{i}\eta)\ud E\ud \eta\ri| \\
   &\le &(p!\pi)^{-1}L \|\phi^{(p+1)}\|_\infty|\nu|(\R)\lf(\int_0^{\eta_{\min}} \eta^{p-1}\ud \eta+\del \int_0^1  \eta^{p-\al}\ud \eta\ri)\\
   &\le &(p!\pi)^{-1}L\|\phi^{(p+1)}\|_\infty|\nu|(\R)\lf(\f{\eta_{\min}^p}{p}+\f{\del  }{p-\al+1}\ri)
   \eeq
\epr

 \beg{cor}\la{2611141}
Let $p\ge 1$, let $\nu_N$ be a sequence of signed measures.   
Suppose that for some constants $C,D,c>0$,   we have   $$\sup  \lf\{|m_{\nu_N}(z)|\ste |\real(z)|\le K,\,  \f{C}{\sqrt{\log N}}\le \im z\le D\ri\}\le C\mre^{-c\sqrt{\log N } }.$$ Let $\phi_N$ be a sequence of smooth compactly supported  functions. Then there is $C'=C'(c,C)$ \st  $$|\nu_N(\phi_N)|\le \f{C'|\nu_N|(\R)L_N\|\phi_N^{(p+1)}\|_{\infty}}{(p+1)!(\log N)^{p/2}},$$ where $L_N$ is the Lebesgue measure of the support of $\phi_N$.
\en{cor}

 \bpr We   apply Theorem \re{Helffer-Sjostrand} with $\eta_{\min}=\f{C}{\sqrt{\log N}}$, $\del=C\mre^{-c\sqrt{\log N }}$, $\al=0$. 
 \epr
 
  \subsubsection{An application of Hadamard's three circles theorem}
 The following use of Hadamard's three circles theorem is due to Kargin, in \cite{KarginPTRF12}. All ideas of the proof of Theorem \re{KHTCT} can be found in \cite{KarginPTRF12}, but as it is not stated clearly, we give a short proof.

  \beg{Th}\la{KHTCT}Let $a>0$. There is  $\del_0=\del_0(a)>0$ \st  for all $\del\in (0, \del_0)$, for all signed measure $\nu$, $$\sup\{|m_\nu(z)|\ste z=\ii a\f{\mre+\mre^{\ii\tta}}{\mre-\mre^{\ii\tta}},\,\tta\in [0,2\pi]\}\le \del\implies \sup_{z\in H_{a,r(\del)}} |m_\nu(z)|\le\mre^{-\sqrt{-c\log \del}},$$ 
  where   $c:={2|\nu|(\R)}/{a}$, \, $
 r(\del):=\mre^{-4\sqrt{-c/\log \del}}$ and 
 for $a,r>0$,  $H_{a,r }$ denotes the disc with diameter  $\lf[\ii a\f{1-r }{1+r }, \ii a\f{1+r }{1-r }\ri]\subset \C^+,$  \ie the disc with center $\ds \ii a\f{1+r^2}{1-r^2}$ and radius $a\f{2r}{1-r^2}$.
 
  \en{Th}
  
 \bpr The starting point of the proof is the so-called \emph{Hadamard  three circles theorem} \cite{hardy}, stating that  for $f:\dD:=\{\xi\in \C\ste |\xi|<1\}\to\C$   analytic, the function  $\ds M(r):=\sup_{|\xi|=r} |f(\xi)|$ is non decreasing and   the function $\tM(s):= \log(M(\mre^s))$ is convexe on  $(-\infty,0)$. If we suppose moreover that $f$ is such that for a certain constant $c$, for all $r<1$, \be\la{2011143WOD}(1-r)M(r)\le c,\ee then there is $\del_0=\del_0(c)>0$ \st   for all $\del\in (0, \del_0)$, we have  $$ M(\mre^{-1})\le \del\implies 
    M(r_c(\del))\le \eps_c(\del),$$ with $r_c(\del):=\mre^{-4\sqrt{-c/\log \del}}$ and   $\eps_c(\del):= \mre^{-\sqrt{-c\log \del}}$.
Indeed, it is equivalent to prove that there is $m_0=m_0(c)<0$ \st    for all $m<m_0$, we have $$  \tM({-1})\le m\implies    \tM(-4\sqrt{c/|m|})\le {-\sqrt{c|m|}},$$ which follows from the convexity of $\tM$ (applied at $-1<-4\sqrt{c/|m|}<-\sqrt{c/|m|}$).

Then one concludes by noticing that for $a,\nu$ as in the statement of the theorem and  $f(\xi):=m_\nu\lf(\ii a\f{1+\xi}{1-\xi}\ri)$, \eqre{2011143WOD} is satisfied for $c:={2|\nu|(\R)}/{a}$
 \epr

  \beg{cor}\la{corH3CTh} Let $\mu_N, \mu$ be \pro measures \st for a certain $a>0$ and a certain $C>0$, we have $$\sup \{ |m_{\mu_N-\mu}(z)|\ste z=\ii a\f{\mre+\mre^{\ii\tta}}{\mre-\mre^{\ii\tta}},\,\tta\in [0,2\pi]\}\le  CN^{-1}.$$ 
  Then for any $K>0$, there is $N_0=N_0(a,C,K)$ \st for $N\ge N_0$, $$\sup\lf\{|m_{\mu_N-\mu}(E+\ii\eta)|\ste E\in [-K,K], \,   \f{16}{\sqrt{a\log(N/C)}}\le \eta\le a\f{\mre-1}{\mre+1}\ri\} \ \le \ \mre^{-8\sqrt{\log(N/C)/a}}.$$
  \en{cor}
  
  \bpr We apply the previous theorem: here, $c=4/a$, $\del=CN^{-1}$, so that $r(\del)=\mre^{-8/\sqrt{a\log(N/C)}}$. Thus $1-r(\del)\le  {8}/{\sqrt{a\log(N/C)}}$ and it is easy to see that for $N$ large enough, the disc $H_{a,r(\del)}$   contains the set in question here. 
  \epr

 \noindent{\bf Acknowledgments}: we would like to thank the referees for their careful readings of our paper and for all the remarks and suggestions.

 \begin{thebibliography}{10}
     \bibitem{agz} G.~Anderson, A.~Guionnet, O.~Zeitouni \emph{An Introduction to Random Matrices}. Cambridge studies in advanced mathematics, {118} (2009).
\bibitem{BaoErdosSchnelli1} Z. Bao, L. Erd{\H o}s, K. Schnelli \emph{Local Stability of the Free Additive Convolution}.  To appear in J. Funct. Anal. (2016), doi:10.1016/j.jfa.2016.04.006.

\bibitem{BaoErdosSchnelli2} Z. Bao, L. Erd{\H o}s, K. Schnelli \emph{Local law of addition of random matrices on optimal scale}.   arXiv:1509.07080. 

       \bibitem{BD13} A. Basak, A. Dembo \emph{Limiting spectral distribution of sums of unitary and orthogonal matrices}. Electron. Commun. Probab. 18 (2013), no. 69, 19 pp.
      \bibitem{BA} {S.T. Belinschi}  \emph{The Lebesgue decomposition of the
free additive convolution of two probability distributions},
Probab. Theory Related Fields, Vol. 142 (2008), no. 1-2, 125--150. 
\bibitem{SerbanLinfty} S.T. Belinschi \emph{$L^\infty$-boundedness of density for free additive convolutions}, Rev. Roumaine Math. Pures
Appl. 59(2), 173--184 (2014).
\bibitem{bbg07} S.T. Belinschi, F. Benaych-Georges, A. Guionnet \emph{Regularization by free additive convolution, square and rectangular cases}. Complex Anal. Oper. Theory 3 (2009), no. 3, 611--660.
\bibitem{BB07} S.T. Belinschi, H. Bercovici \emph{A new approach to subordination results in free probability}. J. Anal. Math. 101 (2007), 357--365.
\bibitem{bbcf12}  S.T. Belinschi, H. Bercovici, M. Capitaine, M. F\'evrier 
\emph{Outliers in the spectrum of large deformed unitarily invariant models}
 arXiv:1207.5443.
\bibitem{BENRECT} F. Benaych-Georges \emph{Rectangular random matrices, related convolution}
Probab. Theory Related Fields, Vol. 144, no. 3 (2009) 471--515. 
\bibitem{BEN} F. Benaych-Georges \emph{Exponential bounds for the support convergence in the Single Ring Theorem}, J. Funct. Anal., Vol. 268 (2015), pp. 3492--3507. 
\bibitem{FBGAK2016} F. Benaych-Georges, A. Knowles \emph{Lectures on the local semicircle law for Wigner matrices}, arXiv.  
\bibitem{BV1} H. Bercovici, D. Voiculescu \emph{Regularity questions for free convolution, Non selfadjoint operator algebras, operator theory, and related topics}, 37--47. Oper. Theory Adv. Appl. 104 (1998)
     \bibitem{BianePFE} P. Biane \emph{Processes with free increments}. Math. Z. 227 (1998), no. 1, 143--174. 
     \bibitem{BOR1} C. Bordenave, D. Chafa\"{\i} \emph{Around the circular law}, Probab. Surv. 9 (2012), 1--89.
     \bibitem{BourgadeCirc} P. Bourgade, H.-T. Yau, J. Yin  \emph{The local circular law for random matrices},    Probab. Theory Related Fields, Vol. 159, no. 3-4, 545--595 (2014).
     \bibitem{BourgadeCirc2} P. Bourgade, H.-T. Yau, J. Yin  \emph{The local circular law II: the edge case},    Probab. Theory Related Fields, Vol. 159, no. 3-4, 619--660 (2014).
   \bibitem{CGLP}  D. Chafa\"\i, O. Gu\'edon, G. Lecu\'e, A. Pajor  \emph{Interactions between compressed sensing random matrices and high dimensional geometry}. Panoramas et Synth\`eses, 37. Soci\'et\'e Math\'ematique de France, Paris, 2012.
    \bibitem{ESRY} L. Erd\H{o}s, J. Ram\'{\i}rez, B. Schlein, H.-T. Yau  \emph{Universality of sine-kernel for Wigner matrices with a small Gaussian perturbation}. 
Electron. J. Probab. 15 (2010), no. 18, 526--603. 
   \bibitem{ESY2} L. Erd\H{o}s, B. Schlein, H.-T. Yau  \emph{Semicircle law on short scales and delocalization of eigenvectors for Wigner random matrices}, Ann. Prob. 37 (2009).

     \bibitem{KantoTheo} O.P. Ferreira, B.F. Svaiter \emph{Kantorovich's Theorem on Newton's Method}, arXiv:1209.5704.
     \bibitem{GUI} A. Guionnet, M. Krishnapur, O. Zeitouni \emph{The Single Ring Theorem}. Ann. of Math. (2) 174 (2011), no. 2, 1189--1217.
\bibitem{GUI2} A. Guionnet, O. Zeitouni \emph{Support convergence in the Single Ring Theorem}. Probab. Theory Related Fields 154 (2012), no. 3-4, 661--675.
\bibitem{haag2}   U. Haagerup, F.  Larsen \emph{Brown's spectral distribution measure for $R$-diagonal elements in finite von Neumann algebras},  {J.. Funct. Anal.}, {176} (2000), 331--367.
\bibitem{hardy} G.H. Hardy \emph{The mean value of the modulus of ananalytic function}. Proceedings of the London
Mathematical Society, 14:269--277, 1915.
\bibitem{KarginPTRF12} V. Kargin \emph{A concentration inequality and a local law for the sum of two random matrices.}  Probab. Theory Related Fields 154 (2012), no. 3-4, 677--702.
\bibitem{KarginAOP13} V. Kargin \emph{An inequality for the distance between densities of free convolutions.} Ann. Probab. 41 (2013), no. 5, 3241--3260.
\bibitem{KarginAOP13Sub} V. Kargin \emph{Subordination of the resolvent for a sum of random matrices}. Ann. Probab., to appear.
\bibitem{ledoux-amsbook}
M.~Ledoux 
\emph{The concentration of measure phenomenon}.
 Providence, RI, AMS, 2001.
 \bibitem{ns06} A. Nica, R. Speicher \emph{Lectures on the combinatorics of free probability}. London Mathematical Society Lecture Note Series, 335. Cambridge University Press, Cambridge, 2006.
 \bibitem{RUD} M. Rudelson, R. Vershynin  \emph{Invertibility of random matrices: unitary and orthogonal perturbations},    J. Amer. Math. Soc.  27 (2014), 293--338.
  \bibitem{TaoVuKCL} T. Tao, V. Vu \emph{Random matrices: universality of ESDs and the circular law}, Ann. Probab. 38 (2010), 2023--2065, with an appendix by Manjunath Krishnapur.
\bibitem{vdn91} D.V. Voiculescu, K. Dykema, A. Nica  \emph{Free random variables} CRM Monograghs Series No.1, Amer. Math. Soc., Providence, RI, 1992 
\bibitem{yinLCL3} I. Yin \emph{The local circular law III: general case}. Probab. Theory Relat. Fields (2014) 160:679--732.
 \en{thebibliography}

\end{document}